\documentclass[a4paper,11pt]{article}
\usepackage{amscd, amsmath, amssymb, amsfonts, epic}
\usepackage{graphicx}
\usepackage{a4wide}
\usepackage{xypic}
\usepackage{paralist}
\usepackage{booktabs}
\usepackage{multirow}

\title{Extremal rational elliptic threefolds}
\author{Arthur Prendergast-Smith}

\date{}

\def\Z{\text{\bf Z}}
\def\Q{\text{\bf Q}}

\def\P{\text{\bf P}}

\def\And{\text{ and }}
\def\Or{\text { or }}
\def\arrow{\rightarrow}

\def\iso{\cong}

\def\Pic{\text{Pic}} 
 
\def\Vert{\text{Vert}}

\renewcommand{\baselinestretch}{1.0}
\renewcommand{\arraystretch}{1.2}

\newtheorem{theorem}{Theorem}[section]
\newtheorem{lemma}[theorem]{Lemma}
\newtheorem{corollary}[theorem]{Corollary}

\begin{document}

\maketitle

An elliptic fibration is a proper morphism $f:X \arrow Y$ of normal
projective varieties whose generic fibre $E$ is a regular curve of
genus $1$. The Mordell--Weil rank of such a fibration is defined to be
the rank of the finitely generated abelian group $\Pic^0 \ E$ of
degree-$0$ line bundles on $E$. In particular, $f$ is called {\it
  extremal} if its Mordell--Weil rank is $0$.

The simplest nontrivial elliptic fibration is a rational elliptic
surface $f: X \arrow \P^1$. There is a complete classification of
extremal rational elliptic surfaces, due to Miranda--Persson in
characteristic $0$ \cite{MirandaPersson1986} and W. Lang in positive
characteristic \cite{Lang1991,Lang1994}. (See also Cossec--Dolgachev
\cite[Section 5.6]{CossecDolgachev1989}.) The purpose of the present
paper is to produce a corresponding classification of a certain
class of extremal rational elliptic threefolds. For reference, the
results are shown in Table \ref{table_intro}.

Let us say a bit more about exactly which objects we are classifying.
It is a classical fact that any rational elliptic surface is the
blowup of $\P^2$ at the base locus (a $0$-dimensional subscheme of
degree $9$) of a pencil of cubic curves. This description allows one
to compute the Mordell--Weil rank in terms of reducibility properties
of curves in the pencil \cite[Theorem 5.2]{Totaro2008}. In dimension
$3$, the analogous situation is to consider a net ($2$-dimensional
linear system) of quadric surfaces in $\P^3$. The base locus of such a
net is a $0$-dimensional subscheme of degree $8$.  We will see below
that, under a certain nondegeneracy assumption on the net, blowing up
at the base locus gives an elliptic fibration $f: X \arrow \P^2$, and
now we can compute the Mordell--Weil rank of $f$ in terms of
reducibility properties of quadrics in the net. To exploit this, we
will consider in this paper only elliptic threefolds obtained by
blowing up the base locus of a net of quadrics in $\P^3$. Table
\ref{table_intro} gives a list of all nets of quadrics (up to
projective equivalence) which give rise to extremal elliptic
threefolds in this way.

The classification may be of interest for several reasons. Firstly, it
is a natural counterpart of the results of Miranda--Persson and Lang
on extremal rational elliptic surfaces. It is perhaps surprising to
see that the situation for threefolds, in which the classification
contains only a small finite number of cases, is simpler than that for
surfaces. Secondly, the method of proof uses the theory of root
systems in an essential way. This gives a further demonstration of
the strong connection --- elaborated in \cite{DolgachevOrtland1988}
and \cite{CossecDolgachev1989} --- between root systems and
configurations of points in projective space. Finally, the
classification provides `test specimens' for the Cone Conjecture in
birational geometry \cite[Conjecture 8.1]{Totaro2008}. That conjecture
predicts that the threefolds appearing in the classification should be
particularly simple from the point of view of birational geometry.
(More precisely, they have finitely generated Cox ring.)
We will not explore this direction in the present work, but we plan to
do so in a forthcoming paper.

The main results of the paper are as follows. Theorem
\ref{thm_countingrank} relates the Mordell--Weil rank of an elliptic
fibration obtained from a net of quadrics to reducibility properties
of quadrics in the net. Theorem \ref{thm_configs_of_fibres} shows
that, for an extremal fibration, the configuration of reducible
quadrics in the net is constrained by a (fixed) finite root system.
These two theorems combine to yield Theorem
\ref{thm_classification_of_reducible_fibres}, which gives a list of
the possible configurations of reducible quadrics for an extremal
fibration. In Section \ref{section_stdforms} we use the combinatorial
data produced by Theorem \ref{thm_classification_of_reducible_fibres}
to determine all extremal nets up to projective equivalence. Finally
in Section \ref{section_extremalquartics} we relate our extremal
elliptic threefolds to extremal quartic plane curves, via the
discriminant.

\vspace{11pt}
Thanks to Burt Totaro for many helpful comments and
suggestions, and to the reviewer for several interesting additions.

\begin{table}
\begin{center}
\begin{tabular}{|c|c|c|c|}
\hline
{\bf Root lattice} & {\bf $\Pic^0(E)$} & {\bf Type of net}
& {\bf Standard Form} 
 \\
\toprule

 & & & $Q_1=Z^2$ 
\\
$E_7$& $0$ & $\{8\}_1$ & $Q_2=X(Y+W)+YW$\\
&&& $Q_3=XZ+(Y+W)^2 $\\ 

\hline

\multirow{6}{*}{$A_7$} &  \multirow{6}{*}{$\Z/2\Z$} & \multirow{3}{*}{$\{8\}_2$}
& $Q_1=YZ+W^2$ \\
&&&  $Q_2=XZ+YW$\\
&&& $Q_3=XW-Y^2+Z^2$\\ \cline{3-4}

&& \multirow{3}{*}{$\{4,4\}_1$}& $Q_1=ZW$\\
&&&$Q_2=XZ+YW$ \\
&&&$Q_3=XY+Z^2+W^2$\\

\hline

\multirow{6}{*}{$D_6 \oplus A_1$} & \multirow{6}{*}{$\Z/2\Z$} &
\multirow{3}{*}{$\{6,2\}$} &  $Q_1=YZ$ \\
&&&$Q_2=XZ+W^2$ \\
&&&$Q_3=XY+Z^2$\\ \cline{3-4}


&&\multirow{3}{*}{$\{4,4\}_2$}& $Q_1=XY$\\
&&&$Q_2=Z^2$\\
&&&$Q_3=(X+Y)Z+W^2$
\\

\hline

\multirow{6}{*}{$A_5 \oplus A_2$} & \multirow{6}{*}{$\Z/3\Z$} &
\multirow{3}{*}{$\{5,3\}$} &  $Q_1=YZ$ \\
&&&$Q_2=XW+Z^2$ \\
&&&$Q_3=XY+W^2$\\ \cline{3-4}

&&\multirow{3}{*}{$\{3,3,2\}_1$}& $Q_1=YZ$\\
&&&$Q_2=X(Z+W)$\\
&&&$Q_3=XY+W^2$
\\

\hline

&&& $Q_1=X(Y+Z)$\\
$D_4 \oplus 3 A_1$& $(\Z/2\Z)^2$ & $\{4,2,2\}$&$Q_2=YZ$\\
&&&$Q_3=(X+Y)Z+W^2$\\

\hline

\multirow{9}{*}{$2A_3 \oplus A_1$} & \multirow{9}{*}{$\Z/4\Z$} &
\multirow{3}{*}{$\{4,4\}_3$} &  $Q_1=XY$\\
&&&$Q_2=XZ+W^2$\\
&&&$Q_3=YW+Z^2$\\ \cline{3-4}

&&\multirow{3}{*}{$\{3,3,2\}_2$} & $Q_1=XY$\\
&&&$Q_2=ZW$\\
&&&$Q_3=(X+Y)Z+W^2$\\ \cline{3-4}

&&\multirow{3}{*}{$\{2,2,2,2\}$}&$Q_1=XY$\\
&&&$Q_2=ZW$\\
&&&$Q_3=(X+Y)(Z+W)$\\

\hline

&&&$Q_1=(X+Y+Z)W$\\
$7A_1$ & $(\Z/2\Z)^3$ & $\{1,1,1,1,1,1,1,1\}$ &$Q_2=(X+Y+W)Z$\\
&& ($\operatorname{char} k=2$ only) &$Q_3=(X+Z+W)Y$\\

\hline
\end{tabular}
\caption{List of extremal nets. The root lattices and Mordell--Weil
  groups are obtained in Section \ref{section_root}. The admissible
  types of nets are obtained in Section \ref{section_comb}. Standard
  forms are obtained in Section \ref{section_stdforms}.}
\label{table_intro}
\end{center}
\end{table}

\normalsize

\vspace{11pt}

{\bf Notation, conventions, definitions:} We work throughout over an
algebraically closed field $k$. In general the characteristic of $k$
is not specified, though in some contexts we will exclude
characteristics $2$ and $3$. 

The term {\it extremal fibration} will always refer to an extremal
elliptic fibration $f: X \arrow \P^2$ obtained by blowing up the base
locus (in the sense described below) of a net of quadrics in $\P^3$
which satisfies Assumption 1 below. A net of quadrics is called {\it
  extremal} if the corresponding morphism $X \arrow \P^2$ is an
extremal fibration. 

If $Q_1$, $Q_2$, $Q_3$ are quadrics in $\P^3$, we write $\left<
Q_1,Q_2,Q_3 \right>$ to denote the net they span: that is, $ \left<
Q_1,Q_2,Q_3 \right> = \{ \lambda_1 Q_1 + \lambda_2 Q_2 + \lambda_3 Q_3
: \lambda_i \in k, \ \lambda_1, \ \lambda_2, \ \lambda_3 \text{ not
  all } 0 \}$. Similarly, $\left<Q_1,Q_2 \right>$ denotes the pencil
spanned by $Q_1$ and $Q_2$.

A {\it basepoint} of a net $N$ of quadrics can refer to either a point
$p \in \P^3$ in the set-theoretic intersection $\cap_{Q \in N} Q$ of
all quadrics in the net, or a common tangent direction of the net (of
any order). If we intend only a point $p \in \cap_{Q \in N} Q$ then we
will use the term \emph{$\P^3$-basepoint}. The multiplicity of a
$\P^3$-basepoint $p_i$ will be denoted by $m_i$. A net $N$ is of {\it
  type} $\{m_1,\ldots,m_n\}$ if it has $\P^3$-basepoints $p_1, \ldots,
p_n$ of multiplicities $m_1, \ldots, m_n$.

We will use the notation $X_{m_1\cdots m_n}$ to denote a threefold
obtained from $\P^3$ by blowing up at the base locus of any extremal
net of type $\{m_1,\ldots,m_n\}$. Note that for a given type
$\{m_1,\ldots,m_n\}$, there may exist nonisomorphic
spaces $X_{m_1\cdots m_n}$.

We abuse terminology by using the term {\it rank-$2$ quadric} to refer
to a quadric in $\P^3$ which is the union of two distinct planes, even
in characteristic $2$.

We denote by $h$ the pullback to $X$ of the hyperplane divisor class
on $\P^3$, and by $e_i$ the pullback to $X$ of the exceptional divisor
$E_i$ of the blowup of the basepoint $p_i$ ($i=1,\ldots,8$).  For
brevity, we will denote the class $h-e_i-e_j-e_k-e_l$ by $h_{ijkl}$,
and the class $e_i-e_j$ by $e_{ij}$ or sometimes (for clarity)
$e_{i,j}$. We denote by $l$ the class in $N_1(X)$ represented by the
pullback of a line in $\P^3$, and by $l_i$ the class of the pullback
of a line in the exceptional divisor $e_i$.

\section{Preliminaries} \label{prelims}

In this section we explain how to obtain an elliptic fibration from a
net of quadrics in $\P^3$, under a certain nondegeneracy assumption on
the net. We then point out some simple consequences of this assumption
which we will use later in the paper.

First let us consider what restriction is needed on a net of
quadrics in $\P^3$ to ensure that it gives an elliptic fibration as
defined above. Given any net with a chosen set of generators, say $N =
\left< Q_1,Q_2,Q_3 \right>$, we get a rational map $\P^3
\dashrightarrow \P^2$: explicitly, the map is $p \mapsto [Q_1(p),
Q_2(p), Q_3(p)]$.  This map is defined outside the base locus of $N$,
so we would like to `blow up at the base locus' (in some sense) to get
a morphism $f: X \arrow \P^2$ from a smooth threefold to $\P^2$.
Furthermore, since we are interested in elliptic fibrations, we want
the generic fibre of $f$ to be a smooth curve of genus $1$. If the
base locus of the net is reduced (that is, it consists of $8$ distinct
points), we can blow up these $8$ points in the usual way, and we do
in fact get an elliptic fibration. But the condition of reduced base
locus is too restrictive for our purposes --- it is proved in
\cite{PrendergastSmith2009} that there is only one such net which
gives an extremal fibration --- so we would like to relax it as far as
possible.

Consider however the net spanned by the following $3$ quadrics in
$\P^3$ with homogeneous coordinates $[X,Y,Z,W]$:

\begin{align*}
Q_1 = X(X-W), \quad Q_2 = Y(Y-W), \quad Q_3 = ZW.
\end{align*}
This net has $4$ basepoints of multiplicity $1$ at
$[X,Y,Z,W]=[0,0,0,1]$, $[1,0,0,1]$, $[0,1,0,1]$, $[1,1,0,1]$, and a
basepoint of multiplicity $4$ at $p=[0,0,1,0]$. Therefore we get a
rational map $\P^3 \dashrightarrow \P^2$ defined outside these $5$
points. We want to resolve the indeterminacy of this rational map to
get a morphism $f: X \arrow \P^2$ which is an elliptic fibration.
Suppose we are in the characteristic $0$ case: then we can blow up
along points and curves to get a morphism (though not uniquely).
Bertini's theorem then tells us that the general fibre of $f$ is
smooth. On the other hand, the general fibre is birational to a
quartic curve $C = Q \cap Q'$, the intersection of $2$ quadrics in the
net. One can check that any such $C$ is singular at $p$, and hence is
rational. Therefore the general fibre of $f$ is rational.
 
Since we are only interested in elliptic fibrations, we want to
exclude troublesome examples like this one. What went wrong? The
problem is that the differentials $dQ_1$ and $dQ_2$ are both zero at
$p$, so no intersection $Q \cap Q'$ of $2$ quadrics in the net can be
smooth at $p$. Since the generic fibre of $f: X \arrow \P^2$ is
birational to a singular quartic of the form $Q \cap Q'$ (a rational
curve), we never get an elliptic fibration in this case. Therefore in
what follows we assume that all nets of quadrics in $\P^3$ satisfy the
following assumption.

\vspace{11pt}

{\bf Assumption 1:} \emph{There exist quadrics $Q$, $Q'$ in the
  net such that the intersection $Q \cap Q'$ is smooth at the base
  locus of the net. Equivalently, for each $\P^3$-basepoint $p$ of the net,
  there is at most one quadric in the net singular at $p$.}

\vspace{11pt} Under this assumption we obtain an elliptic fibration as
follows. Choose a quartic curve of the form $C=Q \cap Q'$ which is
smooth at the base locus, and a quadric $Q''$, not in the pencil
spanned by $Q$ and $Q'$, also smooth at the base locus. (This is
possible since smoothness at a given point is an open condition on
quadrics.) Since $C$ is smooth, its higher tangent directions uniquely
define the basepoints infinitely near to any multiple basepoint of the
net.  Blowing up repeatedly at these basepoints, we obtain a threefold
$X$ on which the proper transforms of $C$ and $Q''$ are disjoint, and
hence a morphism $f: X \arrow \P^2$. 

For $f:X \arrow Y$ the blowup of a point in a smooth variety of
dimension $n$, we have the formula $K_X = f^*(K_Y) + (n-1)E$, where
$E$ is the exceptional divisor of the blowup. Applying this in the
case where $X$ is obtained from $\P^3$ by blowing up 8 points, we get
$K_X=-4h+2e_1+\cdots+2e_8$.  So the class
$-\frac{1}{2}K_X=2h-e_1-\cdots-e_8$ is represented by the proper
transform on $X$ of any quadric in the net smooth at the base locus.
This means that the morphism $f: X \arrow \P^2$ from the previous
paragraph is the same as the one given by the basepoint-free linear
system $|-\frac{1}{2}K_X|$.  The generic fibre $E$ of $f$ need not be
smooth, but it is a regular scheme. Also, adjunction tells us the
canonical bundle $K_E$ is trivial, so $E$ has arithmetic genus $1$.
Hence $f$ is an elliptic fibration, as claimed.

\vspace{11pt}

{\bf Remark:} It is customary to refer to a fibration as above whose
generic fibre is regular but not smooth as a {\it quasi-elliptic
  fibration}, but since the arguments of this paper apply equally well
in both the elliptic and quasi-elliptic cases, we abuse terminology
and refer to both as elliptic fibrations. Many facts about
quasi-elliptic fibrations are known: for instance, they exist only if
the base field has characteristic $2$ or $3$; also, the geometric
generic fibre $E(\overline{k(\P^2)})$ is always a cuspidal rational
curve \cite[Proposition 5.1.2]{CossecDolgachev1989}. Note that the
final net in Table \ref{table_intro}, which is extremal only in
characteristic $2$, gives a quasi-elliptic fibration.

\vspace{11pt}

{\bf Remark:} It is a classical fact that the fibrations $f:X \arrow
\P^2$ correspond to nets of cubic curves in the plane. In one
direction, projecting from one basepoint of our net $N$ of quadrics
transforms the net of quartic curves in $\P^3$ dual to $N$ to a net of
cubic curves in $\P^2$ with $7$ basepoints; in the other, blowing up
the $7$ basepoints of such a net and taking the universal family
$\mathcal{X}$ of elliptic curves over the resulting surface, we get an
elliptic fibration $\mathcal{X} \arrow \P^2$ birational to our original
fibration $f:X \arrow \P^2$. For more details on this correspondence
see \cite[Section 6.3.3]{DolgachevNotes}.

\vspace{11pt}

Here are some straightforward consequences of Assumption 1. 

\begin{lemma} \label{lemma_nondegeneracy} Given a net of quadrics
  satisfying Assumption 1, no $3$ of the basepoints are collinear, nor
  any $5$ coplanar. More precisely, suppose $X$ is the threefold
  obtained from such a net by blowing up its base locus as described
  above.. Then no class $l-\sum_{k=1}^3 l_{i_k}$ in $N_1(X)$ or
  $h-\sum_{k=1}^5 e_{j_k}$ in $N^1(X)$ is represented by an effective
  cycle.
\end{lemma}

{\bf Proof:} For any choice of distinct indices we have $-K_X \circ
(l-\sum_{k=1}^3 l_{i_k}) = -1$ (where $\circ$ denotes intersection of
cycles on $X$), but this is impossible for an
effective cycle since $-K_X$ is
basepoint-free. 

For the second claim, suppose there was an effective cycle
$h-\sum_{k=1}^5 e_{j_k}$ in $N^1(X)$; its image in $\P^3$ would be a
plane $P$. Choose any quartic curve $C=Q \cap Q'$, an
intersection of $2$ quadrics in the net, which is smooth at the base
locus; such a curve exists by Assumption 1. Its proper transform
$\tilde{C}$ on $X$ has class $4l- \sum_{i=1}^8 l_i$.  Therefore
$(h-\sum_{k=1}^5 e_{j_k}) \circ \tilde{C}=-1$, implying that any such
$C$ is contained in $P$. But smoothness of $C$ at
a finite set of points is an open condition on $Q$ and $Q'$, so this
is impossible. \quad (QED Lemma \ref{lemma_nondegeneracy})

\begin{lemma} \label{lemma_reduciblequads} Given a net of quadrics
  satisfying Assumption 1, we have the following facts:

  \begin{itemize}
  \item There is at most $1$ double plane in the net.
  \item There are at most $n$ irreducible cones with vertices at
  basepoints of the net, where $n$ is the number of distinct
  $\P^3$-basepoints of the net.
\item There are finitely many rank-$2$ quadrics in the net.
\end{itemize}

\end{lemma}

{\bf Proof:} Any double plane is singular at all $\P^3$-basepoints, so
by Assumption 1 we get the first claim. For the second, Assumption $1$
implies there is at most $1$ cone with vertex at a given
$\P^3$-basepoint $p_i$. 

For the final claim, suppose there is a curve of rank-$2$ quadrics in
the net. Then every pencil in the net contains a reducible quadric,
hence its base locus is a reducible quartic in $\P^3$. But each fibre
of $f: X \arrow \P^2$ is birational to the base locus of some pencil
in the net, hence must be reducible. This contradicts regularity of
the generic fibre. \quad (QED Lemma
\ref{lemma_reduciblequads})

\section{Rank of the elliptic fibration} \label{sect_rank}

In this section, we derive a formula for the rank of an elliptic
fibration $f: X \arrow \P^2$ obtained from a net of quadrics in
$\P^3$, in terms of the number of distinct $\P^3$-basepoints of the
net and the number of quadrics of rank $2$ in the net. This
generalises \cite[Theorem 7.2]{Totaro2008}, which gives the 
formula for a net with $8$ distinct $\P^3$-basepoints.

\vspace{11pt}

\begin{theorem} \label{thm_countingrank}
Suppose $f:X \rightarrow \P^2$ is an elliptic fibration arising from a
net of quadrics in $\P^3$. Then the rank $\rho$ of the Mordell--Weil
group of the generic fibre of $f$ is given by

\begin{align*}
\rho = n-d-1
\end{align*}
where $n$ is the number of distinct $\P^3$-basepoints of the
net, and $d$  the number of quadrics of rank $2$ in the net. In
particular, $f$ is extremal if and only if $d=n-1$.
\end{theorem}

{\bf Proof:} The rank of an elliptic threefold $f: X \arrow S$ is
given by the Shioda--Tate--Wazir formula \cite[Theorem
2.3]{HulekKloosterman2008}.  Let us derive this formula in our case
$S=\P^2$. To do this, we imitate the proof of \cite[Theorem
7.2]{Totaro2008}. We have a surjective homomorphism $r: \Pic \ X
\arrow \Pic \ E$ given by restriction of divisors, and so
$\operatorname{rank} \Pic \ E = \operatorname{rank} \Pic \ X -
\operatorname{rank} \operatorname{ker} r$. Since we know that $\Pic \
E = \Pic^0 \ E \, \oplus \, \mathbf{Z}$, this gives
$\operatorname{rank} \Pic^0 \ E =\operatorname{rank} \Pic \ X -
\operatorname{rank} \operatorname{ker} r -1 = 8 -\operatorname{rank}
\operatorname{ker} r$. So we need to calculate the rank of the kernel
of the restriction homomorphism.

The kernel of $r$ is generated by the classes of all irreducible
divisors in $X$ which do not map onto $\P^2$ under $f$. If $\lambda$
is the class of a line in $\P^2$, then  $f^*(\lambda)
= -{1 \over 2} K_X$, so the pullback of any irreducible divisor in $\P^2$
is a multiple of $-{1 \over 2}K_X$. Therefore the kernel of the restriction
homomorphism is generated by $-{1 \over 2}K_X$ together with $r_F$
classes for every irreducible divisor $F$ in $\P^2$ whose preimage in
$X$ consists of $r_F+1$ irreducible components, say
$\sum_{j=1}^{r_F+1}m_{F_j}D_{F_j}$. I claim that the
divisors $D_{F_j}$ for any $F$ and $1 \leq j \leq r_F$ are linearly
independent in $\Pic\ X \otimes \Q$. This follows from the
corresponding fact about a morphism from a surface to a curve
\cite[Lemma II.8.2]{Barthetal2003}, by restricting to the inverse
image of a general line in $\P^2$. So the Mordell--Weil group $\Pic^0 \
E$ has rank $8-1-\sum r_F$. We must show this  can be written as
$n-d-1$, where $n$ is the number of distinct $\P^3$-basepoints of the
net, and $d$ the number of rank-$2$ quadrics in the net.

The map $f: X \arrow \P^2$ is given by resolving the indeterminacy of
the rational map $\P^3 \dashrightarrow \P^2 : p \mapsto [Q_1(p),
Q_2(p), Q_3(p)]$ where $Q_i$ is any (fixed) basis for the net of
quadrics. So a fibre of $f$ is (at least away from the base locus of
the net) the intersection $Q \cap Q'$ of $2$ quadrics in the net,
hence a quartic curve. Let us refer to the corresponding quartic curve
$Q \cap Q'$ in $\P^3$ as the {\it pseudofibre} of $f$ over the given
point.

If the intersection $Q \cap Q'$ is smooth at the base locus, then the
pseudofibre $Q \cap Q'$ is isomorphic to the corresponding fibre of
$f$. If such a fibre contains a line, then this must be the line
through $2$ of the basepoints $p_i$. So there are only finitely many
fibres smooth at the base locus which contain a line. The only other
possibility for a reducible pseudofibre smooth at the base locus is
that it be the union $C_1 \cup C_2$ of $2$ smooth conic curves in
$\P^3$. But each curve $C_i$ is contained in a plane $P_i$ in $\P^3$;
the union $P_1 \cup P_2$ is therefore a rank-$2$ quadric in the net
which is smooth at the base locus.

Note this implies in particular that if a reducible divisor $\Delta$
in $\P^3$ contains a pseudofibre smooth at the base locus and maps to
a curve in $\P^2$, then in fact it maps to a line in $\P^2$. To see
this, assume without loss of generality that $\Delta$ is a union of
pseudofibres. Every pseudofibre contained in $\Delta$ and smooth at
the base locus is contained in some rank-$2$ quadric $Q$, whose image
in $\P^2$ is a line, and Lemma \ref{lemma_reduciblequads} shows there
are finitely many such $Q$. These pseudofibres are dense in $\Delta$,
so the image of $\Delta$ is contained in a finite union of lines in
$\P^2$. If different pseudofibres were contained in different rank-$2$
quadrics, the image of $\Delta$ would be a union of distinct lines,
hence reducible, but this contradicts our assumption.  Therefore the
image of $\Delta$ in $\P^2$ is a line, as required.

So the only possibilities for reducible pseudofibres which are smooth
at the base locus are exactly those described in \cite{Totaro2008}.
Let us therefore consider pseudofibres $Q \cap Q'$ which are not
smooth at the base locus.

Suppose that $Q$ is a quadric in the net smooth at the base locus,
and a pseudofibre $Q \cap Q'$ is singular at a $\P^3$-basepoint $p_i$. This
means that the differentials $dQ$ and $dQ'$ are linearly dependent at
$p_i$, so (multiplying by a constant if necessary) $d(Q-Q')=0$ at
$p_i$. By Assumption 1, this implies that $Q-Q'$ is the unique quadric
$Q_i$ in the net singular at $p_i$, or put another way $Q' = \lambda Q
+ \mu Q_i$. So the pseudofibre $Q \cap Q'$ is singular at $p_i$ if and only
if $Q'$ belongs to the pencil $\lambda Q + \mu Q_i$, implying
that $Q \cap Q' = Q \cap Q_i$. 

Now suppose $C \subset \P^2$ is a curve over which all pseudofibres of
$f$ are singular at a $\P^3$-basepoint $p_i$. Fix a quadric $Q$ in the
net which is smooth at the base locus. Over any point of $f(Q) \cap C$
the pseudofibre of $f$ is singular at $p_i$. Over a point $q \in f(Q)
\cap C$ the pseudofibre is an intersection $Q \cap Q'$, and by the
previous paragraph we can take $Q'=Q_i$. Therefore $q = f(Q) \cap
f(Q_i)$. This holds for all $q \in f(Q) \cap C$, so we have $f(Q) \cap
C = f(Q) \cap f(Q_i)$. Since this is true for any quadric $Q$ in the
net smooth at $p_i$ (which $Q$ comprise a Zariski-open set in the
net), we must have $C=f(Q_i)$.  We conclude that the only subvarieties
of $\P^2$ over which all the pseudofibres of $f$ are singular at the
base locus are the lines $f(Q_i)$, the images of the finitely many
quadrics $Q_i$ in the net singular at the base locus.  

Suppose $D$ is a reducible effective divisor in $X$ whose image $f(D)
\subset \P^2$ is an irreducble curve $C$; without loss of generality,
we can assume $D = f^{-1}(f(D))$, that is, $D$ is a union of fibres.
Contracting the exceptional divisors $E_i$ in $X$, the image of $D$ is
an effective divisor $\Delta \in \P^3$. If some pseudofibre contained
in $\Delta$ is smooth at the base locus, then as explained above
$\Delta$ must be a supported on a rank-$2$ quadric in the net. If the
pseudofibre over every point of $C$ is singular at the base locus, the
previous paragraph implies that $C$ must be one of the lines $f(Q_i)$
in $\P^2$, so $\Delta$ is supported on $Q_i$.

We therefore have three types of contribution to the rank of
$\operatorname{ker} r$: first, the class $-{1 \over 2} K_X$; second,
reducible quadrics in the net smooth at the base locus, each of which
adds $1$ to the rank of the kernel; third, the quadrics $Q_i$ singular
at the base locus. Let us anaylse the contribution of these $Q_i$ to
the rank of the kernel.

First, suppose $Q_i$ is an irreducible reduced cone, with vertex at
$p_i$. The corresponding divisor $f^{-1}(f(Q_i))$ on $X$ has $m_i$
components in total, namely the class of the proper transform of the
cone together with $m_i-1$ classes of the form $e_{j,\,j+1}$. The
preimage of any line in $\P^2$ has class $-{1 \over 2}K_X$ in $\Pic \
X$, so the classes of these $m_i$ components sum to $-{1 \over 2}K_X$.
Therefore $Q_i$ contributes $m_i-1$ to the rank of $\operatorname{ker}
r$.

Next suppose that $Q_i$ is a rank-$2$ quadric in the net singular at
the base locus. The singular locus of $Q_i$ is a line in $\P^3$,
therefore contains at most $2$ basepoints of the net by Lemma
\ref{lemma_nondegeneracy}. The corresponding divisor $f^{-1}(f(Q_i))$
on $X$ has $2+(m_i-1)$ components if $Q_i$ is singular at $1$
basepoint $p_i$ and $2+(m_i-1)+(m_j-1)$ components if $Q_i$ is
singular at $2$ basepoints $p_i$ and $p_j$. Again, in both cases the
classes of these components sum to $-{1 \over 2}K_X$. So in the first
case we get a contribution of $1+(m_i-1)$ to the rank of
$\operatorname{ker} r$, and in the second case a contribution of
$1+(m_i-1)+(m_j-1)$.

Finally, consider the case of a non-reduced quadric $Q_i$ --- that is,
a double plane $2P$. In this case, all quadrics in the net except
$Q_i$ must be smooth at the base locus, by Assumption 1. The (reduced)
plane $P$ passes through some subset of the basepoints, including all
of the $\P^3$-basepoints (which are therefore all multiple). The
proper transform of $P$ on $X$ has class $h-e_{i_1}-\cdots-e_{i_j}$ in
$\Pic \ X$, for some set of distinct indices. Therefore the proper
transform of $Q_i$ on X has class $2(h-e_{i_1}-\cdots-e_{i_j})$. On
the other hand, this proper transform must be disjoint from some
smooth fibre $C$, which has class $4l-\sum_i l_i$. We conclude that
$j$, the number of indices in the expression for the class of $P$,
must be $4$.  Again, the divisor $f^{-1}(f(Q_i))$ has class $-{1 \over
  2}K_X= 2h-\sum_i e_i$. We can rewrite this as a sum of effective
classes as follows:

\begin{align*}
2h-\sum_i e_i \quad = \quad 2(h-e_{i_1}-\cdots-e_{i_j}) + \sum_{p_k}
\sum_{p_l} e_{l,\, l+1} + \ R
\end{align*}

where the first sum is taken over the $\P^3$-basepoints $p_k$, the
second over all basepoints $p_l$ infinitely near to $p_k$, except the
highest, and $R$ is a sum of terms of the form $e_{l,\, l+1}$ which
have already appeared in sum. The number of distinct terms in this sum
is $1+ \sum_{p_i} (m_i-1)$, with the sum taken over all
$\P^3$-basepoints $p_i$. Hence the contribution to the rank of
$\operatorname{ker} r$ is $\sum_{\text{all $\P^3$-basepoints} \ p_i}
(m_i-1)$.

(It may help to think about the fibre of $f$ over a general point of
$f(Q_i)$; this is one of the degenerations of elliptic curves
described by Kodaira in \cite{Kodaira1963}. For instance, if our net
has a single basepoint of multiplicity 8 and a double plane $Q_i=2P$,
then the fibre over the generic point of $f(Q_i)$ is a curve of type
III* in Kodaira's notation.)

Let us now show that the above arguments together give the formula
claimed. In the case of no double plane in the net, the total
contribution to the rank of $\operatorname{ker} r$ from quadrics
singular at the base locus is

\begin{align*}
  \sum_{p_i} (m_i-1) + \sum_{p_j} 1+(m_j-1) + \sum_{p_k,\ p_l} 1+
  (m_n-1)+(m_l-1)
\end{align*}

where the first sum is taken over multiple $\P^3$-basepoints at which the
singular quadric is an irreducible cone, the second over multiple
basepoints at which the singular quadric is a rank-$2$ singular at $1$
basepoint, and the third is taken over pairs of multiple $\P^3$-basepoints
both lying on the singular locus of the same rank-$2$ quadric. Since
every multiple $\P^3$-basepoint is of one of these three types, summing we
get

\begin{align*}
 d_{sing} + \sum_{\text{multiple $\P^3$-basepoints} \ p_i} (m_i-1). 
\end{align*}

where $d_{sing}$ is the number of rank-$2$ quadrics in the net
singular at the base locus. Finally including rank-$2$ quadrics smooth
at the base locus, each of which contributes $1$ to the rank, and the
class $-{1 \over 2} K_X$, we get

\begin{align*}
  \operatorname{rank} \operatorname{ker} r = 1+d+\sum_{\text{multiple
      $\P^3$-basepoints} \ p_i} (m_i-1).
\end{align*}

In the case of a double plane in the net, we know that all rank-$2$
quadrics in the net must be smooth at the base locus, so each
contribute $1$ to the rank of $\operatorname{ker} r$, and also that
there are no cones in the net with vertex at a basepoint.  So using
the formula from a few paragraphs back, and including $-{1 \over 2}
K_X$ again, we get $\operatorname{rank} \operatorname{ker} r
=1+d+\sum_{\text{multiple $\P^3$-basepoints} \ p_i} (m_i-1)$. (Recall that in this
case all $\P^3$-basepoints are multiple, so we are summing over the
same set as before.)

Finally computing the rank $\rho$ of $\Pic^0 \ E$ as $ \rho =
8-\operatorname{rank} \operatorname{ker} r$, we get in both cases

\begin{align*}
\rho &= 8- \left( 1+d+\sum_{\text{multiple $\P^3$-basepoints} \ p_i} (m_i-1) \right) \\
&= 7-d-\sum_{\text{all $\P^3$-basepoints} \ p_i} (m_i-1) \\
&= 7-d-8+n \\
&= n-d-1
\end{align*}

as claimed. \quad (QED Theorem \ref{thm_countingrank})

\section{Extremal fibrations and root systems} \label{section_root}

In this section, we will show that the possibilities for an extremal
fibration are constrained by a certain root system. Together with the
rank formula from Section \ref{sect_rank}, this will lead to a
combinatorial classification of extremal fibrations in Section
\ref{section_comb}.

\vspace{11pt}

More precisely, suppose $f: X \arrow \P^2$ is an extremal fibration.
Call an irreducible divisor in $X$ {\it vertical} if it mapped by $f$
to a curve in $\P^2$; we saw in the previous section that the only
vertical divisors are components of divisors $f^{-1}(L)$, where $L$ is
a line in $\P^2$. We will prove that the possible configurations of
vertical divisors are constrained by maximal-rank subsystems of the
root system $E_7$. Before explaining this, let us state the following
lemma. A proof can be found for instance in \cite[Theorem 6.1.2, Table
5]{GorbOnisVin1994}.

\begin{lemma} \label{lemma_e7subsystems}
The only root subsystems of $E_7$ of finite index are the following: 
\begin{inparaenum}[\itshape a\upshape)]
\item $E_7$,
\item $A_7$,
\item $D_6 \oplus A_1$,
\item $A_5 \oplus A_2$,
\item $D_4 \oplus 3A_1$,
\item $2 A_3 \oplus A_1$,
\item $7A_1$. \ $\square$
\end{inparaenum}
\end{lemma}

\vspace{11pt}

We define a bilinear form denoted $\cdot$ on $\Pic \ X $ as follows: 

\begin{align*}
  \Pic \ X \otimes \Pic \ X & \arrow \Z \\
  D_1 \otimes D_2 & \mapsto D_1 \cdot D_2 := D_1 \circ D_2 \circ
  \left( -{1 \over 2} K_X \right)
\end{align*}

where as before $\circ$ means intersection of algebraic cycles on
$X$. For any $D \in \Pic \ X$, we have $D \cdot \left( -{1 \over 2}
  K_X \right)= D \circ (4l-\sum_i l_i)$, so a divisor belongs to the
corank-1 sublattice $K_X^{\perp}$ if and only if it has degree $0$ on
any fibre of $f$. That means the surjection $r : \Pic \ X \arrow
\Pic \ E$ restricts to a surjection $r:K_X^{\perp} \arrow \Pic^0 \ E$.
So the latter group is finite --- that is, $f$ is extremal --- if and
only if the kernel of $r$ has finite index in $K_X^{\perp}$. But the
kernel of $r$ is generated by the classes of vertical divisors. So
given an extremal fibration $X$, the lattice $\Vert(X) \subset \Pic \
X$ spanned by classes of vertical divisors must be a finite-index
sublattice of $K_X^{\perp}$.

\vspace{11pt} 

It is easy to check that the vectors $h_{1234}, \ e_{12}, \ e_{23},
\ldots, e_{78}$ form a system of simple roots of $K_X^{\perp}$ under
the bilinear form defined above, and hence that $K_X^{\perp}$ is
isomorphic to the affine root system $\tilde{E}_7$. At first sight,
the appearance of root systems in this context may seem surprising,
but there is an explanation. The definition above shows that $D_1
\cdot D_2$ actually computes the intersection number of the curves
$D_1 \cap Q \And D_2 \cap Q$ inside $Q$, the proper transform of a
general quadric in the net. Now $f_{| \, Q} : Q \arrow f(Q) \cong
\P^1$ is a rational elliptic surface, so classical results \cite[p.
201]{Barthetal2003} on elliptic surfaces tell us that the intersection
form on the classes of curves lying in fibres of $f_{| \,Q}$ defines
the structure of a root system. Therefore the original form defined on
$\Pic \ X $ also defines a root system. (For an extensive discussion
of the connection between point sets in projective space and root
systems, see \cite[Chapter $5$]{DolgachevOrtland1988}.)

Define the \emph{radical} $\operatorname{Rad} \Lambda$ of a lattice
$\Lambda$ to be the subgroup of elements $\lambda \in \Lambda$ such
that $\lambda \cdot x=0$ for all $x \in \Lambda$. Then
$\operatorname{Rad}(K_X^{\perp})$ is spanned by the class $-{1 \over
  2} K_X$, and $K_X^{\perp}/\operatorname{Rad}(K_X^{\perp}) \ \cong \
\tilde{E}_7 / \operatorname{Rad}(\tilde{E}_7) \ \cong \ E_7$. For any
extremal fibration $X$, the sublattice $\Vert(X) \subset K_X^{\perp}$
spanned by classes of vertical divisors contains the class $-{1
  \over 2}K_X$, so $\Vert(X)/(-{1 \over 2} K_X)$ injects into $E_7$ as
a subsystem of finite index.

Therefore, given any extremal fibration $X$, the root system
$\Vert(X)/(-{1 \over 2} K_X)$ must be one of the $7$ listed in Lemma
\ref{lemma_e7subsystems}. What does this tell us about the possible
configurations of vertical divisors? We noted above that a vertical
divisor in $X$ must map to a line in $\P^2$. Given any line $L \subset
\P^2$, the divisor $f^*(L) \subset X$ has class $-{1 \over 2} K_X$ in
$\Pic \ X$. Suppose that $f^*(L) =-{1 \over 2} K_X = \sum_{i=1}^k m_i
D_i$, with $D_i$ (distinct) irreducible and effective divisors, $m_i$
natural numbers, and $k>1$. The classes $D_i \ (i=1,\ldots k)$ are
linearly independent in $\Pic \ X \otimes \Q$, hence span a
sub-lattice $\Pic \ X$ of rank $k$ which is contained in $\Vert(X)$.
Passing to the quotient $\Vert(X)/(-{1 \over 2} K_X) \subset E_7$, the
images of these classes span a sub-lattice $\Lambda(L)$ of rank $k-1$.
Moreover,by restricting to the preimage of a general line in $\P^2$,
one can check that each such class has $D_i^2=-2$, so in fact their
images span a subsystem.

By connectedness of the fibres of $f$, the Dynkin diagram of
$\Lambda(L)$ is connected. Conversely if $D_1$ and $D_2$ are
components of $f^*(L_1) \And f^*(L_2)$, with the $L_i$ distinct lines
in $\P^2$, we have $D_1 \cdot D_2=0$, because the restrictions of the
$D_i$ to the preimage of a general line in $\P^2$ lie in different
fibres, hence are disjoint. So the connected components $\Gamma_i$ of
the Dynkin diagram of $\Vert(X)/(-{1 \over 2} K_X)$ correspond exactly
to the subsystems spanned by classes of divisors lying over the
finitely many lines $L_i$ in $\P^2$ for which $f^*(L_i)$ is reducible.
Note also that the number of nodes of $\Gamma_i$ is one less that the
number of components of $f^*(L_i)$, since the classes of those
components sum to $-{1 \over 2} K_X \equiv 0$ in $\Vert(X)/(-{1 \over
  2} K_X)$.

The upshot is that to determine the possible configurations
of $f$-vertical divisors in $X$, we need to determine all graphs
obtainable from the Dynkin diagrams of the subsystems in Lemma
\ref{lemma_e7subsystems} by adding $1$ node to each connected
component. There is one extra condition: given a line $L \subset \P^2$
and the corresponding lattice $\Lambda(L) \subset \Vert(X)$ spanned by
classes of irreducible components of $f^*(L)$, we know that
$\Lambda(L)$ is negative semi-definite but not negative definite. (It
contains $-{1 \over 2} K_X$, which has square 0.) Consequently it is
isomorphic to an affine root system of rank $k-1$. So, we must add our
nodes in such a way that each component of the resulting graph is the
Dynkin diagram of some affine root system. (See for instance
\cite{Humphreys1992} for a classification of these.) The result is the
following:

\begin{enumerate}
\item $E_7$: Here we are adding just $1$ node. The only possible
  outcome is $\tilde{E}_7$.
\item $A_7$: Adding $1$ node, we can get either $\tilde{A}_7 \Or
  \tilde{E}_7$.
\item $A_5 \oplus A_2$: For $n \leq 6$, the only allowed way to add a
  node to $A_n$ yields $\tilde{A}_n$. So in this case we get
  $\tilde{A}_5 \oplus \tilde{A}_2$. (Here, and below, the symbol
  $\oplus$ on the right-hand side simply means the disjoint union of
  graphs.)

\item $2A_3 \oplus A_1$: As above, we get $2\tilde{A}_3 \oplus \tilde{A}_1$.

\item $D_6 \oplus A_1$: The only allowed way to add a node to $D_n \
  (n \geq 4)$ yields $\tilde{D}_n$. So here we get
  $\tilde{D}_6 \oplus \tilde{A}_1$.
\item $D_4 \oplus 3A_1$: As above, get $\tilde{D}_4 \oplus 3 \tilde{A}_1$.

\item $7A_1$: As above, get $7 \tilde{A}_1$. 
\end{enumerate}

We can summarise our results as follows.

\begin{theorem} \label{thm_configs_of_fibres} Suppose $f: X \arrow
  \P^2$ is an extremal fibration. Then the lattice
  $\Vert(X)/(-{1 \over 2} K_X)$ is isomorphic to a finite-index
  subsystem of $E_7$. A choice of finite-index subsystem determines
  the configuration of $f$-vertical divisors on $X$, and all
  possibilities are realised.
\end{theorem}

{\bf Proof:} We have already proved the first claim. It remains to
verify the second and third claims.

For the second claim, we must show that the finite-index
subsystem $\Vert(X)/(-{1 \over 2} K_X) \subset E_7$ determines the
configuration of vertical divisors uniquely. In light of the above
discussion, all we need show is that if $\Vert(X)/(-{1 \over 2} K_X)
\iso A_7$, then the configuration of vertical divisors is not
$\tilde{E}_7$. If the configuration were $\tilde{E}_7$, we would have
$\Vert(X)/(-{1 \over 2} K_X) = \tilde{E}_7/(-{1 \over 2} K_X) = E_7$,
contrary to assumption. So the configuration of vertical divisors is
uniquely determined by a choice of subsystem.

The last claim will be verified in Sections \ref{section_comb} and
\ref{section_stdforms}. In Section \ref{section_comb} we will
determine the combinatorial possibilities for a net of quadrics whose
associated configuration of $f$-vertical divisors is a given graph
$\Gamma$ on this list. Then in Section \ref{section_stdforms} we will
exhibit standard forms for each permitted type of net, which shows
in particular that they exist.  \quad (QED Theorem
\ref{thm_configs_of_fibres})

\renewcommand{\arraystretch}{1.62}

\begin{corollary} \label{cor_mordellweil} Suppose that $f: X \arrow
  \P^2$ is an extremal fibration with generic fibre $E$. Then the
  Mordell--Weil group $\Pic^0 \ E$ is determined by the configuration
  of vertical divisors, and is given by the following table. (The
  types corresponding to a given configuration will be derived in
  Section \ref{section_comb}.)

\vspace{11pt}
\begin{center}
\begin{tabular}{lll}
{\bf Vertical divisors} & {\bf $\Pic^0 \ E$} & {\bf Types}   \\
\toprule

$\tilde{E}_7$  & 0 & $\{8\}_1$\\
$\tilde{A}_7$  & $\Z/2\Z$ & $\{8\}_2$, $\{4,4\}_1$\\
$\tilde{D}_6 \oplus \tilde{A}_1$ & $\Z/2\Z$ & $\{6,2\}, \ \{4,4\}_2$ \\
$\tilde{A}_5 \oplus \tilde{A}_2$ & $\Z/3\Z$ & $\{5,3\}, \{3,3,2\}_1$ \\
$2 \tilde{A}_3 \oplus \tilde{A}_1$ & $\Z/4\Z$ & $\{4,4\}_3 , \
\{3,3,2\}_2, \ \{2,2,2,2\}$ \\
$\tilde{D}_4 \oplus 3 \tilde{A}_1$ & $(\Z/2\Z)^2$ & $\{4,2,2\}$  \\
$7 \tilde{A}_1$ & $(\Z/2\Z)^3$ & $\{1,1,1,1,1,1,1,1\}$ \\
\end{tabular}
\end{center}
\end{corollary}

\renewcommand{\arraystretch}{1.2}

{\bf Proof:} We know from the earlier discussion that 

\begin{align*}
  \Pic^0 \ E \ \cong \ \tilde{E}_7 / \Vert(X) \ \cong E_7/\left( \Vert(X)/(-{1
    \over 2} K_X) \right).
\end{align*}

Theorem \ref{thm_configs_of_fibres} shows that the sub-lattice $\Vert(X)/(-{1 \over 2}
K_X)$ is determined by the configuration of vertical
divisors. Moreover, computing the quotients of $E_7$ by its $7$
finite-index sublattices is straightforward, and gives the results shown.
 \quad (QED Corollary \ref{cor_mordellweil})

 \vspace{11pt} 

 Theorem \ref{thm_configs_of_fibres} and Corollary
 \ref{cor_mordellweil} are an analogue of \cite[Theorem
 5.6.2]{CossecDolgachev1989}, which classifies the possible
 configurations of reducible fibres on an extremal rational elliptic
 surface. It is perhaps surprising --- and certainly pleasant --- that
 the result for threefolds is no more complicated than that for
 surfaces.

\section{Combinatorial classification} \label{section_comb}

In this section we use the list of possible configurations of
vertical divisors from Section \ref{section_root} together with
the rank formula of Theorem \ref{thm_countingrank} to determine the
possible types of an extremal net. In fact, the list gives us
more information: given an extremal net with its type and
configuration of vertical divisors, we can say exactly which
classes $D \in \Pic(X)$ are represented by vertical divisors.

\vspace{11pt}

\begin{theorem} \label{thm_classification_of_reducible_fibres} Suppose
  $f: X \arrow \P^2$ is an extremal fibration given by a net $N$ of
  quadrics in $\P^3$. Then the type of $N$ and the classes of
  irreducible vertical divisors in $X$ are (up to permutation of
  indices) one of the cases shown in Table \ref{table_configs}.
\end{theorem}

Note that for some types $\{m_1,\ldots,m_n\}$ we get several possible
configurations of reducible divisors: we use a subscript (as
$\{m_1,\ldots,m_n\}_i$) to distinguish between these. 

\begin{table}
\begin{tabular}{cc}
$\tilde{E}_7$
\begin{picture}(167.5,60)
\put(10,20){\line(1,0){20}}
\put(35,20){\line(1,0){20}}
\put(60,20){\line(1,0){20}}
\put(85,20){\line(1,0){20}}
\put(110,20){\line(1,0){20}}
\put(135,20){\line(1,0){20}}
\put(82.5,22.5){\line(0,1){20}}
\put(7.5,20){\circle{5}}
\put(32.5,20){\circle{5}}
\put(57.5,20){\circle{5}}
\put(82.5,20){\circle{5}}
\put(107.5,20){\circle{5}}
\put(132.5,20){\circle{5}}
\put(157.5,20){\circle{5}}
\put(82.5,45){\circle{5}}

\put(7.5,10){\makebox(0,0){$e_{12}$}}
\put(32.5,10){\makebox(0,0){$e_{23}$}}
\put(57.5,10){\makebox(0,0){$e_{34}$}}
\put(82.5,10){\makebox(0,0){$e_{45}$}}
\put(107.5,10){\makebox(0,0){$e_{56}$}}
\put(132.5,10){\makebox(0,0){$e_{67}$}}
\put(157.5,10){\makebox(0,0){$e_{78}$}}
\put(82.5,55){\makebox(0,0){$h_{1234}$}}

\end{picture}

&
$\tilde{A}_7$
\begin{picture}(167.5,60)

\put(10,20){\line(1,0){20}}
\put(7.5,22.5){\line(0,1){20}}
\put(35,20){\line(1,0){20}}
\put(60,20){\line(1,0){20}}
\put(85,20){\line(1,0){20}}
\put(110,20){\line(1,0){20}}

\put(10,45){\line(1,0){120}}
\put(132.5,22.5){\line(0,1){20}}

\put(7.5,20){\circle{5}}
\put(32.5,20){\circle{5}}
\put(57.5,20){\circle{5}}
\put(82.5,20){\circle{5}}
\put(107.5,20){\circle{5}}
\put(132.5,20){\circle{5}}
\put(7.5,45){\circle{5}}
\put(132.5,45){\circle{5}}

\put(7.5,10){\makebox(0,0){$e_{12}$}}
\put(32.5,10){\makebox(0,0){$e_{23}$}}
\put(57.5,10){\makebox(0,0){$e_{34}$}}
\put(82.5,10){\makebox(0,0){$e_{45}$}}
\put(107.5,10){\makebox(0,0){$e_{56}$}}
\put(132.5,10){\makebox(0,0){$e_{67}$}}
\put(132.5,55){\makebox(0,0){$e_{78}$}}
\put(7.5,55){\makebox(0,0){$c_1$}}

\end{picture}

\\
$\{8\}_1$& $\{8\}_2$
\\

\begin{picture}(167.5,60)

\put(10,20){\line(1,0){20}}
\put(7.5,22.5){\line(0,1){20}}
\put(35,20){\line(1,0){20}}
\put(60,20){\line(1,0){20}}
\put(85,20){\line(1,0){20}}
\put(110,20){\line(1,0){20}}

\put(10,45){\line(1,0){120}}
\put(132.5,22.5){\line(0,1){20}}

\put(7.5,20){\circle{5}}
\put(32.5,20){\circle{5}}
\put(57.5,20){\circle{5}}
\put(82.5,20){\circle{5}}
\put(107.5,20){\circle{5}}
\put(132.5,20){\circle{5}}
\put(7.5,45){\circle{5}}
\put(132.5,45){\circle{5}}

\put(7.5,10){\makebox(0,0){$e_{12}$}}
\put(32.5,10){\makebox(0,0){$e_{23}$}}
\put(57.5,10){\makebox(0,0){$e_{34}$}}
\put(82.5,10){\makebox(0,0){$h_{1235}$}}
\put(107.5,10){\makebox(0,0){$e_{56}$}}
\put(132.5,10){\makebox(0,0){$e_{67}$}}
\put(132.5,55){\makebox(0,0){$e_{78}$}}
\put(7.5,55){\makebox(0,0){$h_{1567}$}}
\end{picture}

&
$\tilde{D}_6 \oplus \tilde{A}_1$
\begin{picture}(167.5,100)

\put(40,45){\line(1,0){40}}
\put(85,45){\line(1,0){40}}

\put(75,69){\line(1,0){20}}
\put(75,71){\line(1,0){20}}

\put(37.5,47.5){\line(0,1){20}}
\put(37.5,42.5){\line(0,-1){20}}
\put(127.5,47.5){\line(0,1){20}}
\put(127.5,42.5){\line(0,-1){20}}

\put(37.5,70){\circle{5}}
\put(37.5,20){\circle{5}}
\put(37.5,45){\circle{5}}
\put(82.5,45){\circle{5}}
\put(127.5,45){\circle{5}}
\put(127.5,70){\circle{5}}
\put(127.5,20){\circle{5}}

\put(72.5,70){\circle{5}}
\put(97.5,70){\circle{5}}

\put(37.5,10){\makebox(0,0){$h_{1278}$}}
\put(27.5,45){\makebox(0,0){$e_{23}$}}
\put(37.5,80){\makebox(0,0){$e_{12}$}}
\put(82.5,35){\makebox(0,0){$e_{34}$}}
\put(127.5,10){\makebox(0,0){$h_{1234}$}}
\put(137.5,45){\makebox(0,0){$e_{45}$}}
\put(127.5,80){\makebox(0,0){$e_{56}$}}
\put(72.5,80){\makebox(0,0){$e_{78}$}}
\put(97.5,80){\makebox(0,0){$c_7$}}

\end{picture}

\\
$\{4,4\}_1$&$\{6,2\}$
\\

\begin{picture}(167.5,100)

\put(40,45){\line(1,0){40}}
\put(85,45){\line(1,0){40}}

\put(75,69){\line(1,0){20}}
\put(75,71){\line(1,0){20}}

\put(37.5,47.5){\line(0,1){20}}
\put(37.5,42.5){\line(0,-1){20}}
\put(127.5,47.5){\line(0,1){20}}
\put(127.5,42.5){\line(0,-1){20}}

\put(37.5,70){\circle{5}}
\put(37.5,20){\circle{5}}
\put(37.5,45){\circle{5}}
\put(82.5,45){\circle{5}}
\put(127.5,45){\circle{5}}
\put(127.5,70){\circle{5}}
\put(127.5,20){\circle{5}}

\put(72.5,70){\circle{5}}
\put(97.5,70){\circle{5}}

\put(37.5,10){\makebox(0,0){$e_{34}$}}
\put(27.5,45){\makebox(0,0){$e_{23}$}}
\put(37.5,80){\makebox(0,0){$e_{12}$}}
\put(82.5,35){\makebox(0,0){$h_{1256}$}}
\put(127.5,10){\makebox(0,0){$e_{78}$}}
\put(137.5,45){\makebox(0,0){$e_{67}$}}
\put(127.5,80){\makebox(0,0){$e_{56}$}}
\put(72.5,80){\makebox(0,0){$h_{1234}$}}
\put(97.5,80){\makebox(0,0){$h_{5678}$}}

\end{picture}
&

$\tilde{A}_5 \oplus \tilde{A}_2$
\begin{picture}(167.5,80)
\put(10,20){\line(1,0){20}}
\put(7.5,22.5){\line(0,1){20}}
\put(35,20){\line(1,0){20}}
\put(60,20){\line(1,0){20}}
\put(82.5,22.5){\line(0,1){20}}
\put(110,20){\line(1,0){20}}
\put(109.5,22){\line(1,1){21}}
\put(132.5,22.5){\line(0,1){20}}

\put(10,45){\line(1,0){70}}

\put(7.5,20){\circle{5}}
\put(32.5,20){\circle{5}}
\put(57.5,20){\circle{5}}
\put(82.5,20){\circle{5}}
\put(107.5,20){\circle{5}}
\put(132.5,20){\circle{5}}

\put(7.5,45){\circle{5}}
\put(82.5,45){\circle{5}}
\put(132.5,45){\circle{5}}

\put(7.5,55){\makebox(0,0){$h_{1678}$}}
\put(7.5,10){\makebox(0,0){$e_{12}$}}
\put(32.5,10){\makebox(0,0){$e_{23}$}}
\put(57.5,10){\makebox(0,0){$e_{34}$}}
\put(82.5,10){\makebox(0,0){$e_{45}$}}
\put(82.5,55){\makebox(0,0){$h_{1234}$}}
\put(107.5,10){\makebox(0,0){$e_{67}$}}
\put(132.5,10){\makebox(0,0){$e_{78}$}}
\put(132.5,55){\makebox(0,0){$c_6$}}

\end{picture}
\\
$\{4,4\}_2$&$\{5,3\}$
\\
\begin{picture}(167.5,80)
\put(10,20){\line(1,0){20}}
\put(7.5,22.5){\line(0,1){20}}
\put(35,20){\line(1,0){20}}
\put(60,20){\line(1,0){20}}
\put(82.5,22.5){\line(0,1){20}}
\put(110,20){\line(1,0){20}}
\put(109.5,22){\line(1,1){21}}
\put(132.5,22.5){\line(0,1){20}}

\put(10,45){\line(1,0){70}}

\put(7.5,20){\circle{5}}
\put(32.5,20){\circle{5}}
\put(57.5,20){\circle{5}}
\put(82.5,20){\circle{5}}
\put(107.5,20){\circle{5}}
\put(132.5,20){\circle{5}}

\put(7.5,45){\circle{5}}
\put(82.5,45){\circle{5}}
\put(132.5,45){\circle{5}}

\put(7.5,55){\makebox(0,0){$h_{1478}$}}
\put(7.5,10){\makebox(0,0){$e_{12}$}}
\put(32.5,10){\makebox(0,0){$e_{23}$}}
\put(57.5,10){\makebox(0,0){$h_{1245}$}}
\put(82.5,10){\makebox(0,0){$e_{56}$}}
\put(82.5,55){\makebox(0,0){$e_{45}$}}
\put(107.5,10){\makebox(0,0){$e_{78}$}}
\put(132.5,10){\makebox(0,0){$h_{1237}$}}
\put(132.5,55){\makebox(0,0){$h_{4567}$}}

\end{picture}

&
$2 \tilde{A}_3 \oplus \tilde{A}_1$
\begin{picture}(167.5,80)
\put(10,20){\line(1,0){20}}
\put(7.5,22.5){\line(0,1){20}}
\put(10,45){\line(1,0){20}}
\put(32.5,22.5){\line(0,1){20}}

\put(60,20){\line(1,0){20}}
\put(57.5,22.5){\line(0,1){20}}
\put(60,45){\line(1,0){20}}
\put(82.5,22.5){\line(0,1){20}}

\put(7.5,20){\circle{5}}
\put(32.5,20){\circle{5}}
\put(7.5,45){\circle{5}}
\put(32.5,45){\circle{5}}

\put(57.5,20){\circle{5}}
\put(82.5,20){\circle{5}}
\put(57.5,45){\circle{5}}
\put(82.5,45){\circle{5}}

\put(107.5,20){\circle{5}}
\put(107.5,45){\circle{5}}

\put(106.5,22.5){\line(0,1){20}}
\put(108.5,22.5){\line(0,1){20}}

\put(7.5,10){\makebox(0,0){$e_{12}$}}
\put(32.5,10){\makebox(0,0){$e_{23}$}}
\put(57.5,10){\makebox(0,0){$e_{56}$}}
\put(82.5,10){\makebox(0,0){$e_{67}$}}
\put(107.5,10){\makebox(0,0){$h_{1234}$}}

\put(7.5,55){\makebox(0,0){$c_1$}}
\put(32.5,55){\makebox(0,0){$e_{34}$}}
\put(57.5,55){\makebox(0,0){$c_5$}}
\put(82.5,55){\makebox(0,0){$e_{78}$}}
\put(107.5,55){\makebox(0,0){$h_{5678}$}}
\end{picture}

\\
$\{3,3,2\}_1$&$\{4,4\}_3$
\\

\begin{picture}(167.5,80)
\put(10,20){\line(1,0){20}}
\put(7.5,22.5){\line(0,1){20}}
\put(10,45){\line(1,0){20}}
\put(32.5,22.5){\line(0,1){20}}

\put(60,20){\line(1,0){20}}
\put(57.5,22.5){\line(0,1){20}}
\put(60,45){\line(1,0){20}}
\put(82.5,22.5){\line(0,1){20}}

\put(7.5,20){\circle{5}}
\put(32.5,20){\circle{5}}
\put(7.5,45){\circle{5}}
\put(32.5,45){\circle{5}}

\put(57.5,20){\circle{5}}
\put(82.5,20){\circle{5}}
\put(57.5,45){\circle{5}}
\put(82.5,45){\circle{5}}

\put(107.5,20){\circle{5}}
\put(107.5,45){\circle{5}}

\put(106.5,22.5){\line(0,1){20}}
\put(108.5,22.5){\line(0,1){20}}

\put(7.5,10){\makebox(0,0){$e_{12}$}}
\put(32.5,10){\makebox(0,0){$e_{23}$}}
\put(57.5,10){\makebox(0,0){$e_{45}$}}
\put(82.5,10){\makebox(0,0){$e_{56}$}}
\put(107.5,10){\makebox(0,0){$e_{78}$}}

\put(7.5,55){\makebox(0,0){$h_{1456}$}}
\put(32.5,55){\makebox(0,0){$h_{1278}$}}
\put(57.5,55){\makebox(0,0){$h_{1234}$}}
\put(82.5,55){\makebox(0,0){$h_{4578}$}}
\put(107.5,55){\makebox(0,0){$c_7$}}
\end{picture}
&
\begin{picture}(167.5,80)
\put(10,20){\line(1,0){20}}
\put(7.5,22.5){\line(0,1){20}}
\put(10,45){\line(1,0){20}}
\put(32.5,22.5){\line(0,1){20}}

\put(60,20){\line(1,0){20}}
\put(57.5,22.5){\line(0,1){20}}
\put(60,45){\line(1,0){20}}
\put(82.5,22.5){\line(0,1){20}}

\put(7.5,20){\circle{5}}
\put(32.5,20){\circle{5}}
\put(7.5,45){\circle{5}}
\put(32.5,45){\circle{5}}

\put(57.5,20){\circle{5}}
\put(82.5,20){\circle{5}}
\put(57.5,45){\circle{5}}
\put(82.5,45){\circle{5}}

\put(107.5,20){\circle{5}}
\put(107.5,45){\circle{5}}

\put(106.5,22.5){\line(0,1){20}}
\put(108.5,22.5){\line(0,1){20}}

\put(7.5,10){\makebox(0,0){$e_{12}$}}
\put(32.5,10){\makebox(0,0){$h_{1356}$}}
\put(57.5,10){\makebox(0,0){$e_{56}$}}
\put(82.5,10){\makebox(0,0){$h_{1257}$}}
\put(107.5,10){\makebox(0,0){$h_{5678}$}}

\put(7.5,55){\makebox(0,0){$h_{1378}$}}
\put(32.5,55){\makebox(0,0){$e_{34}$}}
\put(57.5,55){\makebox(0,0){$h_{3457}$}}
\put(82.5,55){\makebox(0,0){$e_{78}$}}
\put(107.5,55){\makebox(0,0){$h_{1234}$}}
\end{picture}
\\
$\{3,3,2\}_2$&$\{2,2,2,2\}$
\\
$\tilde{D}_4 \oplus 3 \tilde{A}_1$

\begin{picture}(167.5,80)

\put(10,32.5){\line(1,0){20}}
\put(35,32.5){\line(1,0){20}}

\put(32.5,35){\line(0,1){7.5}}
\put(32.5,22.5){\line(0,1){7.5}}

\put(32.5,32.5){\circle{5}}
\put(7.5,32.5){\circle{5}}
\put(57.5,32.5){\circle{5}}
\put(32.5,20){\circle{5}}
\put(32.5,45){\circle{5}}

\put(82.5,20){\circle{5}}
\put(82.5,45){\circle{5}}

\put(81.5,22.5){\line(0,1){20}}
\put(83.5,22.5){\line(0,1){20}}

\put(107.5,20){\circle{5}}
\put(107.5,45){\circle{5}}

\put(106.5,22.5){\line(0,1){20}}
\put(108.5,22.5){\line(0,1){20}}

\put(132.5,20){\circle{5}}
\put(132.5,45){\circle{5}}

\put(131.5,22.5){\line(0,1){20}}
\put(133.5,22.5){\line(0,1){20}}

\put(7.5,22.5){\makebox(0,0){$e_{12}$}}
\put(32.5,10){\makebox(0,0){$h_{1278}$}}
\put(32.5,55){\makebox(0,0){$h_{1256}$}}
\put(42.5,37.5){\makebox(0,0){$e_{23}$}}
\put(57.5,22.5){\makebox(0,0){$e_{34}$}}

\put(82.5,10){\makebox(0,0){$h_{1234}$}}
\put(82.5,55){\makebox(0,0){$h_{5678}$}}

\put(107.5,10){\makebox(0,0){$e_{56}$}}
\put(107.5,55){\makebox(0,0){$c_5$}}

\put(132.5,10){\makebox(0,0){$e_{78}$}}
\put(132.5,55){\makebox(0,0){$c_7$}}

\end{picture}
&
$7 \tilde{A}_1$
\begin{picture}(167.5,80)

\put(7.5,20){\circle{5}}
\put(7.5,45){\circle{5}}

\put(6.5,22.5){\line(0,1){20}}
\put(8.5,22.5){\line(0,1){20}}

\put(32.5,20){\circle{5}}
\put(32.5,45){\circle{5}}

\put(31.5,22.5){\line(0,1){20}}
\put(33.5,22.5){\line(0,1){20}}

\put(57.5,20){\circle{5}}
\put(57.5,45){\circle{5}}

\put(56.5,22.5){\line(0,1){20}}
\put(58.5,22.5){\line(0,1){20}}

\put(82.5,20){\circle{5}}
\put(82.5,45){\circle{5}}

\put(81.5,22.5){\line(0,1){20}}
\put(83.5,22.5){\line(0,1){20}}

\put(107.5,20){\circle{5}}
\put(107.5,45){\circle{5}}

\put(106.5,22.5){\line(0,1){20}}
\put(108.5,22.5){\line(0,1){20}}

\put(132.5,20){\circle{5}}
\put(132.5,45){\circle{5}}

\put(131.5,22.5){\line(0,1){20}}
\put(133.5,22.5){\line(0,1){20}}

\put(157.5,20){\circle{5}}
\put(157.5,45){\circle{5}}

\put(156.5,22.5){\line(0,1){20}}
\put(158.5,22.5){\line(0,1){20}}

\put(7.5,10){\makebox(0,0){$h_{1234}$}}
\put(7.5,55){\makebox(0,0){$h_{5678}$}}

\put(32.5,10){\makebox(0,0){$h_{1256}$}}
\put(32.5,55){\makebox(0,0){$h_{3478}$}}

\put(57.5,10){\makebox(0,0){$h_{1278}$}}
\put(57.5,55){\makebox(0,0){$h_{3456}$}}

\put(82.5,10){\makebox(0,0){$h_{1357}$}}
\put(82.5,55){\makebox(0,0){$h_{2468}$}}

\put(107.5,10){\makebox(0,0){$h_{1368}$}}
\put(107.5,55){\makebox(0,0){$h_{2457}$}}

\put(132.5,10){\makebox(0,0){$h_{1458}$}}
\put(132.5,55){\makebox(0,0){$h_{2367}$}}

\put(157.5,10){\makebox(0,0){$h_{1467}$}}
\put(157.5,55){\makebox(0,0){$h_{2358}$}}

\end{picture}
\\
$\{4,2,2\}$&$\{1,1,1,1,1,1,1,1\}$
\end{tabular}

\caption{Configurations of vertical divisors on extremal
  fibrations.}
 \label{table_configs}
\end{table}
\renewcommand{\arraystretch}{1.0}
\vspace{11pt}

Before proving the theorem, we need some facts about the structure of
the diagram which describes the configuration of vertical divisors.
For brevity, let us denote by $\Gamma_X$ the diagram of irreducible
vertical divisors on an extremal fibration $X$, and by $h^0(\Gamma_X)$
the number of connected components of $\Gamma_X$.

\begin{lemma} \label{lemma_A_1pieces} Suppose $X$ is an extremal
  fibration.  If $\Gamma_X$ has a component $\gamma$ of type
  $\tilde{A}_1$, then the nodes of $\gamma$ are either
\begin{inparaenum}[\itshape a\upshape)]
\item the class $c_i$ of a cone with vertex $p_i$, a basepoint
  of the net with multiplicity $2$, and the class $e_{i,i+1}$, or
\item the classes $h_{abcd}$ and $h_{ijkl}$ of $2$ planes whose union
  is a rank-$2$ quadric in the net, smooth at the base locus (so that 
  $\{\{a,b,c,d\},\{i,j,k,l\}\}$ is a partition of $\{1,\ldots,8\}$).
\end{inparaenum}

In the first case we will say the $\tilde{A}_1$ component is
\emph{conical}; in the second we will say it is \emph{smooth}. 
\end{lemma}

{\bf Proof:} Note that irreducible vertical divisors $D_i$ and $D_j$
are nodes of an $\tilde{A}_1$-component if and only if $D_i \cdot D_j
=2$. To prove the lemma, we simply need to consider the intersection
numbers of different types of vertical divisors.

First consider the class of a double plane. In any net containing a
double plane, all other quadrics are smooth at the base locus (by
Assumption 1). So the only types of vertical divisors are
classes of planes $h_{abcd}$ and  divisors $e_{i,i+1}$.
But given any plane $h_{abcd}$ which is a component of a quadric in
the net, at least one of $\{a,b,c,d\}$ is the index of a
$\P^3$-basepoint. If $h_{ijkl}$ is the class of a double plane, then
all indices of $\P^3$-basepoints are contained in $\{i,j,k,l\}$. So
$\# (\{a,b,c,d\} \cap \{i,j,k,l\}) \geq 1$, hence $h_{abcd} \cdot
h_{ijkl} \leq 1$. Also, $h_{abcd} \cdot e_{i,i+1} = -1, 0 $ or $ 1$ for
any $i$. So the class of a double plane cannot be a node of a
component of type $\tilde{A}_1$.

Next consider $h_{abcd}$, the class of a component of a rank-$2$
quadric $Q$ in the net which is singular at the base locus, say at
$p_a$. Then the other component of $Q$ has class $h_{aijk}$ for some
indices $i, \ j, \ k$, and we have $h_{abcd} \cdot h_{aijk} \leq 1$.
Also, $h_{abcd} \cdot e_{i,i+1} = -1, 0 $ or $ 1$ for any $i$.  Finally, the
class of a singular cone with vertex at $p_i$ is $2h-2e_i-\sum_{k \neq
  i, j} e_k$ (where $p_j$ is the highest-order basepoint infinitely
near to $p_i$).  Calculating, this gives $h_{abcd} \cdot c_i \leq 1$.
So $h_{abcd}$ cannot be a node of component of type
$\tilde{A}_1$.

Next consider the class $c_i$ of a singular cone with vertex at
$p_i$. The argument in the last paragraph shows that
$c_i \cdot h_{abcd} \leq 1$ for any class $h_{abcd}$. If $c_\iota$ is the
class of a cone with vertex at another basepoint $p_\iota$, then $c_i
\cdot c_\iota = (2h-2e_i-\sum_{k \neq i, j} e_k) \cdot (2h-2e_\iota-\sum_{k
  \neq \iota, \lambda} e_k)$. But the first sum includes a term
$e_\iota$ (since $p_\iota$ is not infinitely near to $p_i$) and the
second includes $e_i$. So this is $8-2-2-\#(\{1,\ldots,8\} -
\{i,j,\iota,\lambda\} ) = 0$. Also $c_i \cdot e_{j,j+1} =2$ if and
only if $c_i$ contains a term $-2e_j$, but no term $e_{j+1}$ --- that
is, if and only if $i=j$ and $p_i$ is a basepoint of multiplicity
2. This gives the first case above.

Next consider the class $h_{abcd}$ of a component of a quadric $Q$ in
the net smooth at the base locus. We have seen already that $h_{abcd}
\cdot c_i < 2$ and $h_{abcd} \cdot e_{i,i+1} < 2$ for all $i$. Also,
clearly $h_{abcd} \cdot h_{ijkl}= 2$ if and only if $\{a,b,c,d\} \cap
\{i,j,k,l\} = \emptyset$ --- that is, if and only if $h_{ijkl}$ is the
class of the other component of $Q$.  This gives the second case of
the statement.

Finally, consider a class $e_{i,i+1}$. The only case not yet dealt
with is $e_{i,i+1} \cdot e_{j,j+1}$. Again, one can check that the only
possible values are $-2$, $0$ and $1$. \quad (QED Lemma \ref{lemma_A_1pieces})

\vspace{11pt}

The following lemma was already proved in Section \ref{section_root},
in the discussion preceding Theorem \ref{thm_configs_of_fibres}. We
repeat it here to fix notation, and to emphasise the role of Theorem
\ref{thm_countingrank} in the classification argument which follows.

\begin{lemma} \label{lemma_number_of_components} Let $X$ be an
  extremal fibration. Then the number of components $h^0(\Gamma_X)$ is
  equal to $A+B+C+D$, where

\begin{align*}
A &= \text{number of double planes in the net;} \\
B &= \text{number of rank-$2$ quadrics in the net singular at some
$\P^3$-basepoint;}\\
C &= \text{number of rank-$2$ quadrics in the net smooth at the
base locus;}\\
D &= \text{number of cones in the net with vertex at some $\P^3$-basepoint.}
\end{align*}
In particular, since $B+C=d=n-1$ in the notation of Theorem
\ref{thm_countingrank}, we have $n \leq h^0(\Gamma_X)+1$. \quad $\square$
\end{lemma}

{\bf Proof of Theorem \ref{thm_classification_of_reducible_fibres}:}
We saw in the previous section that the graph $\Gamma_X$ of
irreducible vertical divisors on an extremal fibration $X$ must be one
of the $7$ graphs in Table \ref{table_configs}. To prove the theorem,
we will consider each of these graphs $\Gamma$ in turn, and determine
for which types $\{m_1,\ldots,m_n\}$ of nets there can exist an
extremal net of that type with configuration of vertical divisors
equal to $\Gamma$. This process rests on several earlier results.
First, Theorem \ref{thm_countingrank} tells us how many rank-$2$
quadrics an extremal net of a given type must contain. Next, the first
lemma above narrows down the possibilities for a component of type
$\tilde{A}_1$ in any of the graphs. Finally, the second lemma allows
us to ignore types $\{m_1,\ldots,m_n\}$ with more than $h^0(\Gamma)+1$
distinct $\P^3$-basepoints.

For the purposes of the proof, let us introduce some terminology. A
\emph{simple chain} is a connected graph consisting of nodes
$n_1,\ldots,n_k$, edges (of multiplicity 1) joining $n_i$ to $n_{i+1}$
for $i=1,\ldots k-1$, and no other edges. A \emph{simple $k$-chain} is
a simple chain with $k$ nodes. 

For any net of quadrics in $\P^3$, we adopt the following convention
in labelling its basepoints. First choose a $\P^3$-basepoint, and call
it $p_1$. If $p_1$ has multiplicity $m_1$, then we define $p_2$ to be
the basepoint in the exceptional divisor $E_1$, $p_3$ to be the
basepoint in the exceptional divisor $E_2$, and so on up to $p_{m_1}$.
We then choose $p_{m_1+1}$ to be another $\P^3$-basepoint, and repeat,
until we have exhausted all basepoints. So for instance, if we have a
net of type $\{5,2,1\}$, its $\P^3$-basepoints will be labelled $p_1$,
$p_6$, and $p_8$.

Suppose $Q=P_1 \cup P_2$ is a rank-$2$ quadric in an extremal net,
with $\P^3$-basepoints $p_1,\ldots,p_{i_k}$. We will use the (somewhat
imprecise) notation $Q=1^{m_1}2^{m_2}\ldots
k^{m_n}+1^{\mu_1}2^{\mu_2}\ldots k^{\mu_k}$ to indicate that the plane
$P_1$ (resp. $P_2$) has intersection multiplicities with a smooth
quartic $C=Q_1 \cap Q_2$ ($Q_1, \ Q_2$ quadrics which, together with
$Q$, span the net) equal to $m_1,\ldots,m_n$ (resp.
$\mu_1,\ldots,\mu_k$) at $p_1,\ldots,p_{i_k}$. We refer to such an
expression as the \emph{multiplicity data} of $Q$. Note that there are
various constraints on multiplicity data for rank-$2$ quadrics in the
net. For one, the sums $\sum m_i$ and $\sum \mu_j$ of exponents
appearing in each term must always be $4$, since any plane in $\P^3$
intersects a quartic curve with multiplicity $4$. Also, the
`intersection' of the two terms must consist of at most $2$
basepoints, since if two planes in $\P^3$ share $3$ non-collinear
points $p_i$, they are equal. So for example an expression of the form
$Q=1^22^2+1^22^13^1$ is not permitted.

Now let us consider each graph $\Gamma$ from Table \ref{table_configs}
in turn.

\vspace{11pt}

\begin{enumerate}
\item First consider the case $\Gamma = 7\tilde{A}_1$. I claim that
  the only possible type in this case is $\{1,1,1,1,1,1,1,1\}$. To see
  this, note that the base locus of any net contains at most 4
  multiple basepoints. So at most 4 of the $\tilde{A}_1$-components of
  $\Gamma$ are conical, hence at least $3$ are smooth. So there are at
  least $3$ rank-$2$ quadrics in the net smooth at the base locus. I
  claim any set $\{Q_1,Q_2,Q_3\}$ of $3$ such quadrics must span the
  net.

  If not, the third quadric would belong to the pencil spanned by the
  other $2$; rescaling, we could write $Q_3=Q_1+Q_2$. By assumption
  $Q_1=L_1 \Lambda_1$ and $Q_2=L_2 \Lambda_2$, products of linear
  forms. I claim that the set $\{L_1,\Lambda_1,L_2\}$ is linearly
  independent. If not, we could write $\alpha L_1 + \beta \Lambda_1 +
  \gamma L_2=0$. None of the coefficients in this relation can be
  zero, since by assumption the components of $Q_1$ and $Q_2$ are all
  distinct (they give distinct elements of $\Pic(X)$). So we see that
  $L_1=L_2=0$ implies $\Lambda_1=0$, meaning that $Q_1 \cap Q_2$
  contains a line $L_1=\Lambda_1=0$ along which $Q_1$ is singular.
  Intersecting with any other $Q'$ in the net but not in the pencil
  $\left<Q_1,Q_2 \right>$, we would get a point in the base locus at
  which $Q_1$ is singular, which contradicts the fact that $Q_1$ gives
  an $\tilde{A}_1$-component of smooth type. We conclude that $L_1,
  \Lambda_1, L_2$ are linearly independent. So changing coordinates,
  we can assume that $Q_1=XY$, $Q_2=ZL$, where $L$ is a nonzero linear
  form which is not a multiple of $X$, $Y$ or $Z$. If the coefficient
  of $W$ in $L$ is zero, then both $Q_1$ and $Q_2$ are singular at
  $[0,0,0,1]$, violating Assumption 1. So $L$ must have nonzero
  coefficient of $W$, hence by changing coordinates $W \mapsto L$ we
  get $Q_3=XY+ZW$, which is not reducible. This contradicts our
  assumption, and so we conclude that any such set $\{Q_1, Q_2,Q_3\}$
  must span the net.

  That means that, locally near each basepoint, the base locus $Q_1
  \cap Q_2 \cap Q_3$ of the net is given by the intersection of $3$
  planes. If there were a multiple basepoint $p_i$, then the
  intersection of the $3$ planes at $p_i$ would not be transverse,
  hence not proper. So no multiple basepoint can exist, and the net
  must have type $\{1,1,1,1,1,1,1,1\}$.

\vspace{11pt}

Assume now we have an extremal net of type $\{1,1,1,1,1,1,1,1\}$.
Since there are no multiple basepoints, the $7$ rank-$2$ quadrics in
the net are smooth at the base locus. We must show that the classes of
components of these quadrics are (up to permutation of indices) as
shown in Table \ref{table_configs}.

To see this, note first that there are at most $3$ classes of the form
$h_{12ij}$. If there were $4$ or more, we would have to choose at
least $8$ indices from the set $\{3,4,5,6,7,8\}$. Hence at least one
index would be repeated --- say (by relabelling) the index $3$. Then
there would be $2$ classes of the form $h_{123j}$, which is
impossible. So there at most $3$ classes of the form $h_{12ij}$, and
hence by symmetry at most three classes of the form $h_{abij}$, for
any pair $\{a,b\} \subset \{1,\ldots, 8\}$.

For each rank-$2$ quadric $Q$ in the net, a given basepoint lies in
exactly $1$ component of $Q$, so given an index $a \in \{1,\ldots,
8\}$, exactly $7$ of the $14$ classes $h_{ijkl}$ in the graph have $a
\in \{i,j,k,l\}$. Consider the $7$ classes $h_{aijk}$: there are $21$
indices to choose from $\{1,\ldots,8\} - \{a\}$, with each index
appearing at most $3$ times (by the previous paragraph). The only
possibility is that each index appears exactly $3$ times.

Therefore for any pair $\{a,b\} \in \{1,\ldots, 8\}$, there are
exactly $3$ nodes of the graph which have the form $h_{abij}$.  Since
no $2$ classes $h_{abij}$ can share $3$ indices, each index in the set
$\{1,\ldots,8\} - \{a,b\}$ appears in exactly one of these classes.
Geometrically, this means that given $3$ basepoints $p_a, \ p_b , \
p_c$ of the net, the plane spanned by these $3$ points is a component
of a rank-$2$ quadric in the net, and contains a fourth basepoint
$p_d$ of the net. 

We can relabel basepoints if necessary so that $p_4$
is the fourth basepoint on the plane spanned by $\{p_1,p_2,p_3\}$,
$p_6$ is the fourth basepoint on the plane spanned by $\{p_1, p_2,
p_5\}$, and $p_7$ is the fourth basepoint on the plane spanned by
$\{p_1,p_3, p_5\}$. This gives the classes $h_{1234}$, $h_{1256}$, and
$h_{1357}$ (and since every node in the graph determines the node to
which it is connected, the $3$ classes joined to these) appearing in
the diagram.

To determine the remaining classes, consider the plane spanned by
$\{p_1, p_2, p_7\}$. No two classes $h_{ijkl}$ can share $3$ indices,
so this plane cannot contain $p_3$, $p_4$, $p_5$, or $p_6$.  Therefore
its fourth basepoint must be $p_8$, so there is a node $h_{1278}$.
Similar arguments show we must have nodes $h_{1368}$, $h_{1458}$, and
$h_{1467}$. Since every node in the graph determines the node to which
it is connected, this completes the proof that the nodes of the graph
(possibly after permuting indices) must be the configuration in Table
\ref{table_configs} labelled $\{1,1,1,1,1,1,1,1\}$.
   
\item The next case is $\Gamma=\tilde{D}_4 \oplus 3 \tilde{A}_1$. Here
  $h^0(\Gamma)=4$, so we need only consider types with at most 5
  basepoints. Also, note that if we had a basepoint $p_1$ of
  multiplicity 5 or more, we would have effective divisors
  $e_{12},\ldots,e_{45}$. This would imply that there is a subgraph of
  $\Gamma$ which is a simple $4$-chain. But $\Gamma$ has no such
  subgraph, so we need not consider types with basepoints of
  multiplicity 5 or more.  The remaining types are $\{4,4\}$,
  $\{4,3,1\}$, $\{4,2,2\}$, $\{4,2,1,1\}$, $\{4,1,1,1,1\}$,
  $\{3,3,2\}$, $\{3,3,1,1\}$, $\{3,2,2,1\}$, $\{3,2,1,1,1\}$,
  $\{2,2,2,2\}$, $\{2,2,2,1,1\}$.

\begin{enumerate} 
\item Type $\{4,4\}$: we can rule out this possibility as follows. We
  know $h^0(\Gamma)=A+B+C+D=A+D+n-1=A+D+1$. But $A \leq 1$ and $D \leq
  2$ by Lemma \ref{lemma_reduciblequads}, and $A=1$ implies $D=0$
  (since if there is a double plane in the net, all other quadrics
  must be smooth at the base locus). Therefore $h^0(\Gamma_X) \leq 3$
  for this type of net, so it does not yield $\Gamma$.

\item Type $\{4,3,1\}$: Since this type has no basepoint of
  multiplicity $2$, Lemma \ref{lemma_A_1pieces} says there is no conical
  $\tilde{A}_1$-component. So all $3$ of the $\tilde{A}_1$-components
  are smooth, implying there are at least $3$ rank-$2$ quadrics in the
  net. This is impossible by Theorem \ref{thm_countingrank}, so this
  type does not give $\Gamma$.

\item Type $\{4,2,2\}$: The nodes $e_{12}, \ e_{23}, \ e_{34}$ form a
  simple $3$-chain, which must be contained in the
  $\tilde{D}_4$-component. The nodes $e_{56} $ and $ e_{78}$ are
  disjoint from this chain, and from each other, so they must belong
  to $2$ distinct $\tilde{A}_1$-components, which are therefore
  conical, with nodes $e_{56}, \ c_5$ and $e_{78}, \ c_7$.  Since a
  conical $\tilde{A}_1$-component comes from a basepoint of
  multiplicity exactly 2, the third such component must be smooth. So
  there must be a rank-$2$ quadric in the net smooth at the base
  locus.  Clearly, the only possibility for the multiplicity data is
  $Q=1^4+2^23^2$. The corresponding nodes of the diagram are $h_{1234}
  $ and $ h_{5678}$. The other rank-$2$ quadric in the net must
  therefore be singular at the base locus, and its components must
  give the other $2$ nodes in the $\tilde{D}_4$-component. Suppose a
  class $h_{abcd}$ has $h_{abcd} \cdot e_{12}=0, \ h_{abcd} \cdot
  e_{23} = 1, \ h_{abcd} \cdot e_{34}=0$. Then the set $\{a,b,c,d\}$
  contains 1 and 2, but not 3 or 4. Also, we have $h_{abcd} \cdot
  e_{56} = h_{abcd} \cdot e_{78}=0$, so $\{a,b,c,d\}$ must intersect
  both $\{5,6\}$ and $\{7,8\}$ in either 0 or 2 elements. The only two
  possibilities are $h_{1256}$ and $h_{1278}$. This gives the
  configuration in Table \ref{table_configs} labelled $\{4,2,2\}$.

\item Type $\{4,2,1,1\}$: Here there is only $1$ basepoint of
  multiplicity $2$, so at most $1$ conical $\tilde{A}_1$-component. It
  is easy to see that the only possibility for a smooth
  $\tilde{A}_1$-component is $Q=1^4+2^23^14^1$, so we cannot obtain
  the remaining $2$ such components. Hence a net of this type cannot
  yield $\Gamma$.

\item Type $\{4,1,1,1,1\}$: Here there are no basepoints of
  multiplicity 2, so all the $\tilde{A}_1$-components must be
  smooth. But again the only possibility is $Q=1^4+2^13^14^15^1$, so
  we cannot get $\Gamma$ from a net of this type.

\item Types $\{3,3,2\}$ and $\{3,3,1,1\}$: These  types have $2$
  \emph{disjoint} simple $2$-chains, with nodes $e_{12}, \ e_{23}$ and
  $e_{34}, \ e_{45}$. But there is no way to embed $2$ such chains
  disjointly in $\Gamma$, so these types cannot yield $\Gamma$.

\item Type $\{3,2,2,1\}$: In this case there are $2$ basepoints of
  multiplicity 2, giving nodes $e_{45}$ and $e_{67}$ which are
  disjoint from the $2$-chain with nodes $e_{12}, \ e_{23}$ and from
  each other. So these must give $2$ distinct conical
  $\tilde{A}_1$-components.  Also, there are $3$ rank-$2$ quadrics
  in the net, giving 6 more nodes. Adding all these up, this gives 12
  nodes in total, whereas $\Gamma$ has only 11 nodes. So this type
  cannot yield $\Gamma$.

\item Type $\{3,2,1,1,1\}$: This is similar to the previous case. The
  basepoint of multiplicity 3 gives a simple $2$-chain with nodes
  $e_{12}$, $e_{23}$, which must be contained in the
  $\tilde{D}_4$-component. The basepoint of multiplicity 2 gives a
  node $e_{45}$ disjoint from this, so it must be a node of a conical
  $\tilde{A}_1$-component. The net has 4 rank-$2$ quadrics, giving 8
  more nodes. Adding these up we get 12 nodes, so this type cannot
  give $\Gamma$.

\item Type $\{2,2,2,2\}$: We cannot have 3 smooth
  $\tilde{A}_1$-components, for the same reason as in the case
  $\Gamma=7 \tilde{A}_1$. So one of these components must be conical;
  without loss of generality, we have a cone $c_1$. Then all other
  quadrics in the net are smooth at $p_1$. In particular, the 3
  rank-$2$ quadrics in the net are all smooth at $p_1$. But then
  exactly the same argument as in the case $\Gamma=7 \tilde{A}_1$
  shows the intersection is not proper.

\item Type $\{2,2,2,1,1\}$: Just as in the previous case we must have
  a conical  $\tilde{A}_1$-component, so the 4 rank-$2$
  quadrics in the net must all be smooth at $p_1$ say. Again this
  implies the intersection is not proper.
\end{enumerate}

\item The next graph to consider is $\Gamma = 2\tilde{A}_3 \oplus
  \tilde{A}_1$. It has $h^0(\Gamma)=3$, so we need only consider nets
  with at most $4$ basepoints. A basepoint of multiplicity at
  least $5$ would give a simple $4$-chain embedded in $\Gamma$, so we
  know that all basepoints have multiplicity at most $4$. The remaining
  types are $\{4,4\}$, $\{4,3,1\}$, $\{4,2,2\}$, $\{4,2,1,1\}$,
  $\{3,3,2\}$, $\{3,3,1,1\}$, $\{3,2,2,1\}$, $\{2,2,2,2\}$. 

  \begin{enumerate} 
  
  \item Type $\{4,4\}$: There is no basepoint of multiplicity 2, hence
    the $\tilde{A}_1$-component must be smooth. Its nodes are
    therefore $h_{1234} $ and $ h_{5678}$. Also there are 2 simple
    $3$-chains, with nodes $e_{12}, \ e_{23}, \ e_{34} $ and $ e_{56},
    \ e_{67}, \ e_{78}$. If the net had a double plane, it would have
    class $h_{1256}$. This node would be joined to $e_{23} $ and $
    e_{67}$, giving a component of $\Gamma$ with at least 7 nodes, so
    such a double plane cannot exist. Since the unique rank-$2$
    quadric in the net is smooth at the base locus, we must have cones
    $c_1$ and $c_5$, and these give all nodes of $\Gamma$. The
    resulting diagram is shown in Table \ref{table_configs}, labelled
    $\{4,4\}_3$.

  \item Type $\{4,3,1\}$. Again the $\tilde{A}_1$-component must be
    smooth. Since there is only 1 such component, the other reducible
    quadric in the net is singular at some basepoint. If it were
    smooth at $p_1$, its multiplicity data would be $Q=1^4+2^23^1$,
    hence would be smooth at the base locus. This is impossible, so
    $Q$ must be singular at $p_1$. If it were smooth at $p_5$, its
    multiplicity data would be $Q=1^33^1+1^12^3$, so the corresponding
    nodes would be $h_{1238}$ and $h_{1567}$. The first node would be
    joined to $e_{34}$ and the second to $e_{12}$. This would give a
    component of $\Gamma$ with at least 5 nodes, which does not exist.
    Finally, if $Q$ were singular at both $p_1$ and $p_5$, it would
    look like $Q=1^32^1+1^12^23^1$, giving nodes $h_{1235} $ and $
    h_{1568}$.  But the first node would be joined to both $e_{34}$
    and $e_{56}$, giving a component with at least 6 nodes, again
    impossible. So this type does not yield $\Gamma$.

  \item Type $\{4,2,2\}$: We have a simple $3$-chain with nodes
    $e_{12}$, $e_{23}$, $e_{34}$, which must be contained in one of
    the $\tilde{A}_3$-components, and $2$ other nodes $e_{56}, \
    e_{78}$ disjoint from this $3$-chain and each other. If we assume
    first that the $\tilde{A}_1$-component is conical, then the
    remaining 4 nodes are the components of the 2 rank-$2$ quadrics in
    the net. So we must have a class $h_{abcd}$ having intersection 1
    with both $e_{12} $ and $ e_{34}$ and intersection $0$ with
    $e_{23}$, which is impossible. So we can assume the
    $\tilde{A}_1$-component is smooth; hence its nodes are $h_{1234} $
    and $ h_{5678}$. There are 3 more nodes; $2$ are components of the
    other rank-$2$ quadric in the net. If the third was the class of a
    double plane, it would be $h_{1257}$, which would have
    intersection $1$ with $e_{23}$, which would therefore have degree
    $3$. Since the simple $3$-chain is contained in an
    $\tilde{A}_3$-component, and all nodes of that component have
    degree 2, that is impossible. So the last node must be the class
    of a cone, hence $c_5$ or $c_7$. But either choice would yield
    another double edge of the graph, which is impossible. So this
    type cannot yield $\Gamma$.

  \item Type $\{4,2,1,1\}$: First assume the $\tilde{A}_1$-component
    is conical. Then there are no rank-$2$ quadrics in the net smooth
    at the base locus, hence all 3 rank-$2$ quadrics in the net are
    singular at some basepoint. No quadric in the net is singular at a
    basepoint of multiplicity 1, so each of the $3$ rank-$2$ quadrics
    is singular at one of the 2 multiple basepoints. But then $2$
    quadrics must be singular at the same basepoint, which contravenes
    Assumption 1. So this type cannot yield $\Gamma$.
    
  \item Type $\{3,3,2\}$: Here we have $2$ disjoint simple $2$-chains,
    with nodes $e_{12}$, $e_{23} $ and $ e_{45}$, $e_{56}$, and a node
    $e_{78}$ not joined to either. It follows that $e_{78}$ must be a
    node of the $\tilde{A}_1$-component, which is therefore conical
    with second node $c_7$. The 4 remaining nodes are the components
    of the 2 rank-$2$ quadrics in the net, so each is of type
    $h_{abcd}$. Consider the node $h_{abcd}$ of this type joined to
    $e_{12}$: it is not joined to $e_{23}$, so the set
    $\{a,b,c,d\}$ contains 1 but not 2 or 3.  Also it is disjoint
    from $e_{45} $ and $ e_{56}$, so $\{a,b,c,d\}$ either contains or is
    disjoint from $\{4,5,6\}$. But it cannot be disjoint from a set
    $5$ elements, so we must have $\{a,b,c,d\}=\{1,4,5,6\}$. Similar
    arguments show that the remaining node in this component must be
    $h_{1278}$, and the $2$ missing nodes in the other component are
    $h_{1234} $ and $ h_{4578}$. This gives the configuration shown in
    Table \ref{table_configs} labelled $\{3,3,2\}_2$.

  \item Type $\{3,3,1,1\}$: There is no basepoint of multiplicity 2,
    so the $\tilde{A}_1$-component must be smooth. So its nodes
    (possibly after swapping $p_7$ and $p_8$) are $h_{1237}$ and
    $h_{4568}$. We have $2$ simple $2$-chains with nodes $e_{12}, \
    e_{23} $ and $ e_{45}, \ e_{56}$; the 4 remaining nodes must be the
    components of the 2 remaining rank-$2$ quadrics in the net. If a
    class $h_{abcd}$ is joined to $e_{12}$ but not $e_{23}$, then 1
    belongs to $\{a,b,c,d\}$, but 2 and 3 do not. Also, $h_{abcd}$ is
    not connected to $e_{45}$ or $e_{56}$, so either $\{4,5,6\}
    \subset \{a,b,c,d\}$, or they are disjoint. They cannot be
    disjoint since 2 and 3 are not in $\{a,b,c,d\}$ also; therefore
    the node connected to $e_{12}$ is $h_{1456}$. But then $h_{1456}
    \cdot h_{1237} =1$, which gives an illegal edge of the graph. So
    this type does not yield $\Gamma$.

  \item Type $\{3,2,2,1\}$: Again we have a simple $2$-chain with
    nodes $e_{12}, \ e_{23}$, and $2$ nodes $e_{45}, \ e_{67}$ not
    connected to that chain or each other. If the
    $\tilde{A}_1$-component was conical, we would have 5 nodes. The
    components of the 3 rank-$2$ quadrics in the net would give a
    further 6, making 11 altogether, which is a contradiction. So the
    $\tilde{A}_1$-component must be smooth. The only possibility for
    the multiplicity data of the corresponding rank-$2$ quadric is
    $Q=1^34^1+2^23^2$, so the nodes of this $\tilde{A}_1$-component
    must be $h_{1238} $ and $ h_{4567}$.

    The $2$-chain must be contained in an $\tilde{A}_3$-component, so
    there must be a node joined to $e_{12}$ but not $e_{23}$. Since
    $p_1$ has multiplicity 3, the class of a cone in the net with
    vertex at $p_1$ would be $c_1=2h-2e_1-e_2-e_4-\cdots -e_8$, which
    would give $c_1 \cdot e_{23}=1$. So the node in question must be
    the class of a plane $h_{abcd}$. By the same logic as before, the
    set $\{a,b,c,d\}$ contains 1 but not 2 or 3, and since $h_{abcd}
    \cdot e_{45} = h_{abcd} \cdot e_{67} =0$, it must contain or be
    disjoint from the sets $\{4,5\} $ and $ \{6,7\}$. So after
    relabelling basepoints if necessary, it is $h_{1458}$. One can
    check that the final node in that component of $\Gamma$ must be
    $h_{1267}$.

    There are $2$ remaining nodes with classes $h_{abcd} $ and $
    h_{ijkl}$, which must both be joined to $e_{45} $ and $ e_{67}$, but
    no other nodes. So $\{a,b,c,d\} $ and $ \{i,j,k,l\}$ both contain 4
    and 6 but not 5 or 7. Since neither node is joined to $e_{12}$ or
    $e_{23}$, the $2$ sets must also contain or be disjoint from
    $\{1,2,3\}$. But neither is possible. So this type does not yield
    $\Gamma$.

  \item Type $\{2,2,2,2\}$: As before, if we have a conical
    $\tilde{A}_1$-component, say with node $c_1$, then the 3 rank-$2$
    quadrics in the net are smooth at $c_1$, so the intersection is
    not proper. So the $\tilde{A}_1$-component must be smooth.
    Possibly relabelling basepoints, its nodes are $h_{1234} $ and $
    h_{5678}$. The 4 components of the remaining rank-$2$ quadrics in
    the net (which must be singular at the base locus) give the 4
    remaining nodes. Suppose $e_{12}$ and $e_{56}$ belonged to the
    same $\tilde{A}_3$-component. Then there would be a node
    $h_{abcd}$ joined to $e_{12}$ and $e_{56}$, and to no other nodes.
    So the set $\{a,b,c,d\}$ must contain 1 and 5, but not 2 or 6;
    also, it must either contain or be disjoint from $\{3,4\}$ and
    $\{7,8\}$. If it contained $\{3,4\}$, then the intersection
    $h_{abcd} \cdot h_{1234}$ would be $-1$; if it was disjoint from
    $\{3,4\}$, the intersection would be $1$. Since $h_{1234}$ belongs
    to the $\tilde{A}_1$-component, neither is possible. The same
    argument shows that $e_{12}$ and $e_{78}$ cannot belong to the
    same component. So, one $\tilde{A}_3$-component contains $e_{12}$
    and $e_{34}$, and the other contains $e_{56}$ and $e_{78}$.

    So there are $2$ nodes $h_{abcd}$ joined to $e_{12}$ and $e_{34}$
    and no other nodes: it is not hard to see they must be $h_{1356}$
    and $h_{1378}$. Similarly there are $2$ nodes joined to $e_{56}$
    and $e_{78}$ and no other nodes: they must be $h_{1257}$ and
    $h_{3457}$.  This gives the configuration shown in Table
    \ref{table_configs} labelled $\{2,2,2,2\}$. 

  \end{enumerate}

\item The next graph to consider is $\Gamma = \tilde{A}_5 \oplus
  \tilde{A}_2$. This has $h^0(\Gamma)=2$, so we need only consider
  types with at most 3 basepoints. Also, the maximum length of a
  simple chain embedded in this graph is 5, so there can be no
  basepoint of multiplicity more than 6. Also note that if a net has a
  basepoint $p_i$ of multiplicity at least $5$, then all rank-$2$
  quadrics in the net must be singular at that basepoint (otherwise we
  would have a smooth quadric intersecting a plane with multiplicity
  at least $5$ at $p_i$). If the net has $3$ basepoints, it has $2$
  rank-$2$ quadrics, which therefore must both be singular at $p_i$.
  But this contradicts Assumption 1. So we can ignore the types
  satisfying these two conditions, namely $\{6,1,1\}$ and $\{5,2,1\}$.
  This leaves the following types to be considered: $\{6,2\}$,
  $\{6,1,1\}$, $\{5,3\}$, $\{5,2,1\}$, $\{4,4\}$, $\{4,3,1\}$,
  $\{4,2,2\}$, $\{3,3,2\}$.

  \begin{enumerate}

  \item Type $\{6,2\}$: The unique rank-$2$ quadric in the net must
    have multiplicity data $Q=1^4+1^22^2$, so it is singular at $p_1$
    and smooth at $p_7$. Because it is singular at $p_1$ there is no
    double plane in this net; because it is smooth at $p_7$, there
    must be a cone in the net with vertex at $p_7$. But this would
    give a node joined to $e_{78}$ by a double edge, which $\Gamma$
    does not possess. So this type does not yield $\Gamma$.

  \item Type $\{5,3\}$: We have a simple $4$-chain with nodes
    $e_{12},\ldots,e_{45}$, and a simple $2$-chain with nodes $e_{67},
    \ e_{78}$, disjoint from it. The longer chain must be contained in
    the $\tilde{A}_5$-component, and the shorter one in the
    $\tilde{A}_2$-component. The third node of the
    $\tilde{A}_2$-component cannot be a class $h_{abcd}$, since we
    cannot have $h_{abcd} \cdot e_{67}= h_{abcd} \cdot e_{78}=1$. (If
    we did, we would have $h_{abcd} \cdot (e_6-e_8) =2$, which is
    impossible.) So it must be the class of a cone, hence $c_6$. The
    $2$ remaining nodes of the $\tilde{A}_5$-component must be the
    components of the unique rank-$2$ quadric $Q$ in the net. The
    multiplicity data must be $Q=1^4+1^13^3$, so these nodes are
    $h_{1234}$ and $h_{1678}$.  This gives the configuration shown in
    Table \ref{table_configs} labelled $\{5,3\}$.

  \item Type $\{4,4\}$: We know that
    $h^0(\Gamma_X)=A+B+C+D=A+D+(n-1)=A+D+1$, so to get
    $h^0=2$ we need $A+D=1$. First suppose $A=0$. There is a unique
    rank-$2$ quadric $Q$ in the net; the only possibilities for the
    multiplicity data are $Q=1^4+2^4 $ and $ Q=1^32^1+1^12^3$. So $Q$
    is singular at neither or both of the basepoints. If neither, then
    there must be cones in the net with vertices at both basepoints,
    hence $D=2$. If both, then there are no cones singular at the base
    locus, hence $D=0$. Neither case gives $h^0=2$. On the other hand,
    if $A=1$, the unique reducible reduced quadric in the net must be
    smooth at the base locus, so $\Gamma_X$ must have an
    $\tilde{A}_1$-component, which $\Gamma$ does not possess. So this
    type does not yield $\Gamma$.

  \item Type $\{4,3,1\}$: We have a simple $3$-chain with nodes
    $e_{12}, \ e_{23}, \ e_{34}$, which must be contained in the
    $\tilde{A}_5$-component. So the simple $2$-chain with nodes
    $e_{56}, \ e_{67}$ must be contained in the
    $\tilde{A}_2$-component. There are 4 more nodes, which are
    therefore the components $h_{abcd}$ of the 2 rank-$2$ quadrics in
    the net. One of these must be the last node of the
    $\tilde{A}_2$-component, so must have $h_{abcd} \cdot e_{56} =
    h_{abcd} \cdot e_{67}=1$. As above this is impossible, so this
    type does not yield $\Gamma$.

  \item Type $\{4,2,2\}$: We have a simple $3$-chain with nodes
    $e_{12}$, $e_{23}$, $e_{34}$, and $2$ nodes $e_{56}$, $e_{78}$ not
    joined to this chain or each other. Possibly after relabelling,
    $e_{78}$ is a node of the $\tilde{A}_2$ component. Again counting
    nodes, the remaining 2 nodes of this component must be classes
    $h_{abcd}$ and $h_{ijkl}$ of components of rank-$2$ quadrics in
    the net. As before, $\{a,b,c,d\} $ and $ \{i,j,k,l\}$ must both
    contain or be disjoint from $\{1,2,3,4\}$. So these index sets
    must be $\{1,2,3,4\}$ and $\{5,6,7,8\}$, and therefore $h_{abcd}
    \cdot h_{ijkl} =2$, which is impossible. So this type does not
    yield $\Gamma$.

    \item Type $\{3,3,2\}$: Here we have $2$ simple $2$-chains, with
      nodes $e_{12}$, $e_{23}$ and $e_{45}$, $e_{56}$, and a node
      $e_{78}$ disjoint from these chains. The 4 remaining nodes must be
      classes $h_{abcd}$ of components of rank-$2$ quadrics
      in the net. 

      Suppose first that $e_{78}$ is a node of the
      $\tilde{A}_2$-component. The other $2$ nodes of that component
      must be classes $h_{abcd}$ and $h_{ijkl}$, where $\{a,b,c,d\}$
      and $\{i,j,k,l\}$ must both contain or be disjoint from
      $\{1,2,3\}$ and $\{4,5,6\}$, and must both contain 7 but not 8.
      So they are $h_{1237}, \ h_{4567}$. What of the other $2$ nodes?
      One is connected to $e_{12}$ but not $e_{23}$: its class is
      $h_{abcd}$, where the index set contains 1 but not 2 or 3. If
      this node is joined to $e_{45}$ then the index set contains 4
      but not 5 or 6, so the class must be $h_{1478}$. The other node
      is then $h_{1245}$. This gives the configuration shown in Table
      \ref{table_configs} labelled $\{3,3,2\}_1$.

      If the node connected to $e_{12}$ is also connected to $e_{56}$,
      the index set contains 4 and 5 but not 6, and neither or both of
      7 and 8. But this is impossible, since we know it contains $1$
      but not $2$ or $3$. 

      Next suppose that $e_{78}$ belongs to the
      $\tilde{A}_5$-component. There is no way to embed $2$ simple
      $2$-chains and $1$ other node disjointly in this component, so
      in this case one of the $2$-chains must belong to the
      $\tilde{A}_2$-component. Say it is the chain with nodes $e_{45},
      \ e_{56}$: then there is a class $h_{abcd}$ with $h_{abcd} \cdot
      e_{45}= h_{abcd} \cdot e_{56}=1$, which again is impossible, and
      similarly for the other $2$-chain.
    \end{enumerate}

  \item The next graph to consider is $\tilde{D}_6 \oplus
    \tilde{A}_1$. Again we need only consider types with no more than
    3 basepoints, and multiplicities no greater than 6. Also, as
    before we can ignore types $\{6,1,1\}$ and $\{5,2,1\}$. The
    remaining types are $\{6,2\}$, $\{5,3\}$, $\{4,4\}$, $\{4,3,1\}$,
    $\{4,2,2\}$, $\{3,3,2\}$.

    \begin{enumerate}

    \item Type $\{6,2\}$: We have a simple $5$-chain with nodes
      $e_{12}, \ldots, e_{56}$, and a node $e_{78}$ disjoint from it.
      A simple $5$-chain in $\tilde{D}_6$ is joined to all nodes of
      that component, so $e_{78}$ must be a node of the
      $\tilde{A}_1$-component, which is therefore conical, with the
      other node equal to $c_7$. The multiplicity data of the unique
      rank-$2$ quadric in the net must be $Q=1^4+1^22^2$, so
      the corresponding nodes are $h_{1234}$ and $h_{1278}$. This
      gives the configuration shown in Table \ref{table_configs}
      labelled $\{6,2\}$. 

    \item Type $\{5,3\}$: Here there is a simple $4$-chain with nodes
      $e_{12},\ldots,e_{45}$ and a disjoint simple $2$-chain with
      nodes $e_{67}, \ e_{78}$. But it is impossible to embed these
      $2$ chains disjointly in $\Gamma$. So this type does not yield
      $\Gamma$.

    \item Type $\{4,4\}$: We saw before that this type has
      $h^0(\Gamma)=2$ only if there is a double plane in the net. Such
      a plane has class $h_{1256}$. We have nodes $e_{12}, \ e_{23}, \
      e_{34}$, $e_{56}, \ e_{67}, \ e_{78}$; together with $h_{1256}$
      these form the $\tilde{D}_6$-component. The unique rank-$2$
      quadric in the net is smooth, hence gives the
      $\tilde{A}_1$-component with nodes $h_{1234}$ and $h_{5678}$. So
      this type yields the diagram shown in Table \ref{table_configs}
      labelled $\{4,4\}_2$.

    \item Type $\{4,3,1\}$: This type has no basepoint of multiplicity
      2, so the $\tilde{A}_1$-component must be smooth, with nodes
      $h_{1234}$ and $h_{5678}$. We have a simple $3$-chain with nodes
      $e_{12}$, $e_{23}$, $e_{34}$ and a disjoint simple $2$-chain
      with nodes $e_{56}$, $e_{67}$. The $2$ remaining nodes must be
      the classes $h_{abcd}$ of the components of the second rank-$2$
      quadric in the of net. There is a unique way (up to graph
      isomorphism) to embed a simple $3$-chain and a simple $2$-chain
      disjointly in $\tilde{D}_6$, so for $\{i,j\}=\{5,6\}$ or
      $\{7,8\}$, one of the classes $h_{abcd}$ must satisfy $ h_{abcd}
      \cdot e_{ij} =1, \ h_{abcd} \cdot D=0$ for all other nodes $D$
      of the graph. In either case $\{a,b,c,d\}$ contains exactly $2$
      of $\{5,6,7\}$. Also it cannot contain $\{1,2,3,4\}$, so must
      be disjoint from it. But then $\{a,b,c,d\}$ contains at most $3$
      elements, a contradiction. So this type cannot yield $\Gamma$.

    \item Type $\{4,2,2\}$: This graph has a single
      $\tilde{A}_1$-component, so there is some rank-$2$ quadric in
      the net singular at the base locus (and therefore by Assumption
      1 no double plane). The multiplicity data of such a rank-$2$
      quadric has the form $Q=1^i2^j3^k+1^{4-i}2^{2-j}3^{2-k}$.  From
      this we see that $Q$ cannot be singular at both basepoints of
      multiplicity $2$, for if it were we would have
      $Q=1^22^13^1+1^22^13^1$, which cannot happen. So there must be
      cones in the net with vertices at $p_5$ and $p_7$, hence nodes
      $c_5$ and $c_7$ in $\Gamma$. These nodes are joined to $e_{56}$
      and $e_{78}$ respectively by double edges; since there is only
      $1$ double edge in $\Gamma$, we get a contradiction. So this
      type cannot yield $\Gamma$.

    \item Type $\{3,3,2\}$: We have $2$ disjoint simple $2$-chains
      with nodes $e_{12}, \ e_{23}$ and $e_{45}, \ e_{56}$, and a node
      $e_{78}$ disjoint from both chains. Any union of $2$ disjoint
      simple $2$-chains in $\tilde{D}_6$ is joined to every node, so
      we conclude that the node $e_{78}$ must belong to the
      $\tilde{A}_1$-component. This component is then conical with
      other node $c_7$, and then together with the components of
      the $2$ rank-$2$ quadrics in the net we get 10 nodes rather than
      9. So this type does not yield $\Gamma$.

      \end{enumerate}
    
    \item The next graph to consider is $\tilde{A}_7$. Here we need
      only consider types with at most $2$ basepoints, hence the
      possible types are $\{8\}$, $\{7,1\}$, $\{6,2\}$, $\{5,3\}$,
      $\{4,4\}$.
      
      \begin{enumerate}

      \item Type $\{8\}$: We have 7 nodes $e_{12}, \ldots, e_{78}$.
        There are no rank-$2$ quadrics in the net, so the only issue
        is whether the quadric singular at the basepoint is a double
        plane or a cone. If it were a double plane, it would have
        class $h_{1234}$, meaning the node $e_{45}$ in the graph would
        have degree 3. The graph $\tilde{A}_7$ has no such node, so
        the final node must be a cone $c_1$. So the only possibility
        is the configuration shown in Table \ref{table_configs}
        labelled $\{8\}_2$.

      \item Type $\{7,1\}$: The unique rank-$2$ quadric must have
        multiplicity data $Q=1^4+1^32^1$, so the corresponding nodes
        must be $h_{1234}$ and $h_{1238}$. These have intersection
        $-1$, which is impossible. So this type cannot yield
        $\Gamma$ --- indeed, it cannot occur at all. 

      \item Type $\{6,2\}$: Here the unique rank-$2$ quadric has
        multiplicity data $Q=1^4+1^22^2$, so the corresponding nodes
        are $h_{1234}$ and $h_{1278}$. But $h_{1278} \cdot e_{23}=1$,
        so $e_{23}$ has degree 3, which again is impossible for this
        graph. So this type does not yield $\Gamma$.

      \item Type $\{5,3\}$: This type has $h^0(\Gamma_X)=
        A+B+C+D=A+D+(n-1)=A+D+1$. In this case $h^0(\tilde{A}_7)=1$,
        so we must have $A=D=0$. Therefore the unique rank-$2$ quadric
        in the net must be singular at both basepoints. But the only
        possible multiplicity data is $Q=1^4+1^12^3$, so $Q$ is smooth
        at one basepoint. Therefore this type does not yield $\Gamma$,
        or indeed any graph with $h^0 = 1$.

      \item Type $\{4,4\}$: This type has $h^0(\Gamma_X)=
        A+B+C+D=A+D+(n-1)=A+D+1$. In this case $h^0(\tilde{A}_7)=1$,
        so we must have $A=D=0$. Therefore the unique rank-$2$ quadric
        in the net must be singular at both basepoints. The
        multiplicity data of this quadric must be $Q=1^32^1+1^12^3$,
        so the corresponding nodes are $h_{1235}$ and $h_{1567}$.
        Together with the $3$-chains $e_{12}, \ e_{23}, \ e_{34}$ and
        $e_{56}, \ e_{67}, \ e_{78}$, these give the configuration shown
        in Table \ref{table_configs} labelled $\{4,4\}_1$. Note that
        this argument shows in fact that any net of type $\{4,4\}$
        with $h^0(\Gamma_X)=1$ must have $\Gamma_X=\tilde{A}_7$.
        \end{enumerate}

      \item The final graph to consider is $\tilde{E}_7$. Here we need
        only consider types with at most $2$ basepoints. We saw above
        that the type $\{7,1\}$ cannot occur, that the type $\{5,3\}$
        cannot give $h^0(\Gamma_X)=1$, and that a net of type
        $\{4,4\}$ with $h^0(\Gamma_X)=1$ must have $\Gamma_X =
        \tilde{A}_7$. So the only types we need to consider are
        $\{8\}$ and $\{6,2\}$.

        \begin{enumerate}

        \item Type $\{8\}$: As for the case $\Gamma=\tilde{A}_7$, the
          only issue is whether the quadric singular at the basepoint
          is a cone or a double plane. We saw that a cone gives
          $\Gamma_X = \tilde{A}_7$, so it must be a double plane, with
          class $h_{1234}$. So the configuration is as shown in Table
          \ref{table_configs} labelled $\{8\}_1$.

        \item Type $\{6,2\}$: We have a simple $5$-chain with nodes
          $e_{12}, \ldots, e_{56}$. The unique rank-$2$ quadric in the
          net must have multiplicity data $Q=1^4+1^22^2$, so the
          corresponding nodes are $h_{1234}$ and $h_{1278}$. But then
          the nodes $e_{23}$ and $e_{45}$ both have degree 3, which is
          impossible in $\tilde{E}_7$. So this type does not yield
          $\Gamma$. \quad (QED Theorem
          \ref{thm_classification_of_reducible_fibres})

          \end{enumerate}

\end{enumerate}

\section{Standard forms for extremal nets} \label{section_stdforms}

The aim of this section is to find standard forms for extremal nets of
the possible types $\{m_1,\ldots,m_n\}$ determined in Theorem
\ref{thm_classification_of_reducible_fibres}. More precisely, for each
possible configuration $\{m_1,\ldots,m_n\}_i$ of irreducible vertical
divisors shown in Table \ref{table_configs}, we give a unique standard
form for extremal nets whose associated configuration is
$\{m_1,\ldots,m_n\}_i$.

\vspace{11pt}

{\bf Note on characteristic:} We must note at this point that some of
the arguments used to obtain the standard forms below are not valid in
characteristics $2$ and $3$. Therefore, \emph{we claim only that the
  standard forms below exist for nets in $\P^3_k$, where
  $\operatorname{char } k = 0 $ or $ p \geq 5$}. On the other hand, it
is straightforward to check that in each case (except the last) the
given net has the configuration of vertical divisors claimed, and that
the net satisfies Assumption 1, for all characteristics. So our
standard forms prove the existence of extremal nets with each possible
configuration except $\{1,1,1,1,1,1,1,1\}$, in all characteristics.

\begin{enumerate}

\item $\{8\}_1$: the standard form is $Q_1=Z^2, \ Q_2 =X(Y+W)+YW, \
  Q_3=XZ+(Y+W)^2$.

\item $\{8\}_2$: the standard form is $Q_1=YZ+W^2, \ Q_2=XZ+YW, \
  Q_3=XW-Y^2+Z^2$.

\item $\{6,2\}$: the standard form is $Q_1=YZ, \ Q_2=XZ+W^2$,
  $Q_3=XY+Z^2$.

  \item $\{5,3\}$: the standard form is $Q_1=YZ, \ Q_2=XW+Z^2, \ Q_3=XY+W^2$.
     
  \item $\{4,4\}_1$: the standard form is $Q_1=ZW, \ Q_2=XZ+YW, \
    Q_3=XY+Z^2+W^2$.
  \item $\{4,4\}_2$: the standard form is $Q_1=XY, \ Q_2=Z^2, \ Q_3 =
    (X+Y)Z+W^2$.
  \item $\{4,4\}_3$: the standard form is $Q_1=XY, \ Q_2=XZ+W^2, \ Q_3=YW+Z^2$.

  \item $\{4,2,2\}$: the standard form is $Q_1=X(Y+Z), \ Q_2=YZ, \
    Q_3=XZ+W^2$.

 \item $\{3,3,2\}_1$: the standard form is $Q_1=XY, \ Q_2=ZW, \
    Q_3=(X+Y)Z+W^2$.
 \item $\{3,3,2\}_2$: the standard form is $Q_1=YZ, \ Q_2=X(Z+W), \
    Q_3=XY+W^2$.
 
  \item $\{2,2,2,2\}$: the standard form is $Q_1=XY, \ Q_2=ZW, \ Q_3 =
    (X+Y)(Z+W)$.

  \item $\{1,1,1,1,1,1,1,1\}$: extremal nets of this type exist only
    in characteristic $2$, and have standard form $Q_1=(X+Y+Z)W$,
    $Q_2=(X+Y+W)Z$, $Q_3=(X+Z+W)Y$.

\end{enumerate}

The remainder of this section gives a detailed derivation of the
standard forms presented above. 

\begin{enumerate}

\item $\{8\}_1$: First suppose we have a net of type $\{8\}$ which
  contains a double plane. I claim we can put it in standard form
  $Q_1=Z^2, \ Q_2 =XY+XW+YW, \ Q_3=XZ+(Y+W)^2$. To see this, first
  apply a projective transformation moving the unique basepoint to
  $[X,Y,Z,W]=[1,0,0,0]$.  Next, applying an element of $PGL(3) \subset
  PGL(4)$ fixing $p_1$, we can move the double plane so that
  (set-theoretically) it becomes the plane $\{Z=0\}$. This gives $Q_1$
  the form we claimed.

  Next consider $Q_2$. I claim that we can choose $Q_2$ to be an
  irreducible reduced cone with vertex not lying on $Q_1$. To see
  this, consider the subset $S \subset \P^3$ consisting of all singular
  points of all quadrics in the net. I claim $S$ is not contained
  in $\{Z=0\}$.

  Suppose it were, and assume first that the set of singular quadrics
  spans the net. Choose $2$ singular quadrics $Q, \ Q'$ which,
  together with $Q_1$, span the net. By assumption, $Q$ and $Q'$ are
  singular at some point of $\{Z=0\}$. The intersection $Q_1 \cap Q
  \cap Q'$ is a single 8-fold point $p_1$, which means that both $Q_1
  \cap Q$ and $Q_1 \cap Q'$ must be quadruple lines, meeting at $p_1$.
  It is then not difficult to see that we can find a quadric in the
  pencil spanned by $Q$ and $Q'$ which is singular at $p_1$. But this
  violates Assumption 1 above.

  On the other hand, suppose that the set of singular quadrics is
  contained in a pencil. This pencil is spanned by $Q_1$ and any other
  singular quadric $Q_2$, which by assumption is a cone with vertex
  lying in the plane $\{Z=0\}$. We can move the vertex to $p_2 =
  [0,1,0,0]$ without changing $Q_1$ or $p_1$. Adding a multiple of
  $Q_1$ to $Q_2$ does not change the differential at a point of
  $\{Z=0\}$, so every quadric in the pencil (hence every singular
  quadric in the net) is singular at $p_2$ (and nowhere else, unless
  it is $Q=Q_1$). Now choose any smooth quadric $Q_3$ in the net, and
  consider the intersection $Q_1 \cap Q_2 \cap Q_3$.  (Note that $Q_3$
  is not contained in the pencil $\left< Q_1,Q_2 \right>$, so this
  intersection is the set of $\P^3$-basepoints of the net.)
  If $Q_1 \cap Q_2$ consisted (as a set) of $2$ distinct lines $L_1
  \cup L_2$ in $\{Z=0\}$, with $L_1=\{Z=W=0\}$ the line through $p_1$
  and $p_2$, then $L_2 \cap Q_3$ would give another basepoint of the
  net, contradicting our assumption. So $Q_1 \cap Q_2$ must be a
  double line $\{Z=W=0\}$.  Therefore the form defining $Q_2$ looks
  like $Q_2=\alpha XZ+ \beta Z^2 + \gamma ZW + \delta W^2$.
  Subtracting a multiple of $Q_1$, we can assume $\beta=0$; since
  $Q_2$ is an irreducible cone, neither $\alpha$ nor $\delta$ are
  zero. We can therefore scale $X$ and $W$ to obtain
  $\alpha=\delta=1$, without changing $Q_1$, $p_1$, or $p_2$. Now the
  restriction of the form $Q_3$ to the line $\{Z=W=0\}$ must have a
  double root at $p_1$, so the coefficient of $XY$ in $Q_3$ must be
  zero. The coefficient of $X^2$ in $Q_3$ is also zero, since $Q_3$
  passes through $p_1$.  Finally, we can subtract multiples of $Q_1$
  and $Q_2$ from $Q_3$ to make the coefficients of $XZ$ and $Z^2$
  zero, without changing anything.  Note also that since $p_2$ is not
  a basepoint of the net, the coefficient $\epsilon$ of $Y^2$ in $Q_3$
  is nonzero. But now computing the determinant of a general member of
  the net $\lambda_1Q_1+\lambda_2Q_2+\lambda_3Q_3$, we see that the
  discriminant locus is defined by a degree-4 polynomial in the
  $\lambda_i$ which is different from $\lambda_3^4$ --- specifically,
  the coefficient of $\lambda_2^3 \lambda_3$ equals $-\epsilon$, which
  is nonzero. In other words, the set of singular quadrics in the net
  is not contained in the pencil $\{\lambda_3=0\}=\left< Q_1,Q_2
  \right>$, which contradicts our assumption.

  So without loss of generality, we can choose $Q_2$ to be an
  irreducible reduced cone, with vertex not lying on $Q_1$.  Applying
  a projective transformation fixing $p_1$ and $Q_1$, we can bring
  this vertex to the point $p_2=[0,0,1,0]$. This implies that in the
  equation defining $Q_2$, each monomial containing $Z$ has
  coefficient zero. $Q_2$ passes through $p_1$, so the coefficient of
  $X^2$ is zero also: hence $Q_2=b_2XY+c_2XW+d_2Y^2+e_2YW+f_2W^2$, for
  some coefficients $b_2,\ldots,f_2$.
      
  Next we can change coordinates in the plane $\{Z=0\}$, without
  affecting $p_1, \ p_2 \text{ or } Q_1$. So, choose any $2$ points in
  $Q_1 \cap Q_2$ which do not span a line through $p_1$: we can move
  these to $[0,1,0,0] \text{ and } [0,0,0,1]$. In these coordinates
  $d_2=f_2=0$, so we have $Q_2=b_2XY+c_2XW+e_2YW$. Now, $Q_2$ is an
  irreducible reduced cone with vertex not lying on $\{Z=0\}$, so its
  intersection with this plane must be a smooth conic. So none of
  $b_2, \ c_2, \ e_2$ can be zero. Dividing by $b_2$, we can write
  $Q_2=XY+c_2XW+e_2YW$. Changing coordinates $W \mapsto c_2W$, we get
  $Q_2=XY+XW+e_2YW$. Finally, changing coordinates $Y \mapsto e_2^{-1}
  Y, \ W \mapsto e_2^{-1} W$, we get $Q_2=e_2^{-1}(XY+XW+YW)$. (None of
  these coordinate changes affect $p_1, \ p_2, \ Q_1$ or the $2$
  points fixed above.) So we have $Q_2=XY+XW+YW$, as
  claimed.

  Finally we must deal with $Q_3$. First suppose it is a general
  quadric in the net: it has the form
  $Q_3=b_3XY+c_3XZ+d_3XW+e_3Y^2+f_3YZ+g_3YW+h_3Z^2+i_3ZW+j_3W^2$.  We
  know that the plane curves $Q_2 \cap \{Z=0\} \text{ and } Q_3 \cap
  \{Z=0\}$ must have an intersection point of multiplicity 4 at $p_1$.
  This means the following. Suppose we restrict to the affine chart
  $\{X=1\}$ inside $\{Z=0\}$. Then on $Q_2$, we can express $Y$ say as
  a power series in $W$: $Y=p(W)$. Then, substituting $Y=p(W)$ into
  the equation for $Q_3$, we get a power series $q_3(W)$, and the
  condition is that this vanish to order 4 at $W=0$. Since we already
  know $Q_3$ vanishes at $p_1$, his gives $3$ additional equations in
  the coefficients of $Q_3$, namely $d_3=b_3, \ g_3=b_3+2e_3, \
  j_3=e_3$. Applying these conditions, and replacing $Q_3$ by
  $Q_3-b_3Q_2-h_3Q_1$, we can assume $Q_3=e_3W^2 + 2e_3WY + e_3Y^2 +
  i_3WZ + c_3XZ + f_3YZ$, which simplifies to $Q_3=e_3(Y+W)^2 +
  Z(c_3X+ f_3Y+i_3W)$. From this we see that $e_3$ cannot be zero, for
  then $Q_3$ would be reducible. So dividing across, we can assume
  $e_3=1$. Moreover, we see that the differential $dQ_3$ at $p_1$ is
  just $c \, dZ$, so by Assumption 1 we have $c \neq 0$. So, changing
  coordinates $Z \mapsto cZ$, we can assume $c=1$. So we get
  $Q_3=(Y+W)^2 + Z(X+ f_3Y+i_3W)$. Finally, for a given choice of
  coefficients $f_3, \ i_3$, it is straightforward to find a
  projective transformation taking the net spanned by $Q_1, \ Q_2
  \text{ and this } Q_3$ to that spanned by the standard quadrics
  above.
   
\item $\{8\}_2$: The next case is a net of type $\{8\}$ which does not
  contain a double plane. The unique quadric $Q_1$ in the net which is
  singular at $p_1$ is then a reduced irreducible cone. Again putting
  $p_1=[1,0,0,0]$, we see that $Q_1$ is the cone over a smooth conic
  in $\{X=0\} \iso \P^2$. Standard arguments about smooth quadrics
  show that we can change coordinates so that $Q_1=YZ+W^2$.

  Since $p_1$ is a multiple basepoint, there must be some tangent line
  $L \subset T_{p_1} \P^3$ shared by all quadrics in the net. Let
  $p_2$ be the unique point of $Q_1 \cap \{X=0\}$ such that
  $\overline{p_1p_2}$ has tangent direction $L$ at $p_1$. We can apply
  a projective transformation in $PGL(3) \subset PGL(4)$ to bring
  $p_2$ to the point $[0,1,0,0]$.
        
  Now by Assumption 1, for any choice of generators $Q_2, \ Q_3$ of
  the net, the differentials $d \, Q_2$ and $d \, Q_3$ are linearly
  independent at $p_1$. So given a plane $P \subset T_{p_1} \P^3$
  containing the line $L$, we can find a quadric $Q$ in the net with
  $P$ as its tangent plane at $p_1$. In particular we can choose $Q$
  so that its embedded tangent plane at $p_1$ intersects $Q_1$ in a
  double line. Moreover, this property is unchanged if we replace $Q$
  by $Q+\lambda Q_1$ (for any $\lambda \in k)$. So without loss of
  generality, we can assume that $Q_2$ has embedded tangent plane
  intersecting $Q_1$ in a double line $2L$, and that the coefficient
  of $W^2$ in $Q_2$ is zero. This means that $Q_2$ is given by a form
  $Q_2=XZ+ b_2Y^2+c_2YZ+d_2YW+e_2Z^2+f_2ZW$. (We know the coefficient
  of $XZ$ is nonzero, since $Q_2$ is smooth at $p_1$, so we can divide
  across by that coefficient.)

  Now consider $Q_3$. We know that it passes through $p_1$, and that
  its tangent space at $p_1$ contains the line $L$. This implies that
  the coefficients of the monomials $X^2$ and $XY$ in $Q_3$ vanish.
  Also, subtracting appropriate multiples of $Q_1$ and $Q_2$, we can
  assume that the coefficients of $W^2$ and $XZ$ in $Q_3$ also vanish.
  Finally, since $Q_3$ is smooth at $p_1$, the coefficient of $XW$
  must be nonzero, so we can assume it is 1. Putting these facts
  together, we obtain $Q_3=XW
  +e_3Y^2+f_3YZ+g_3YW+h_3Z^2+i_3ZW+j_3W^2$.

  Now we can use the power-series method explained in the previous
  case to obtain equations in the coefficients of $Q_2$ and $Q_3$.
  These are as follows: $b_2=0, \ b_3=0, \ e_3=-d_2, \
  g_3=c_2+c_3d_2-g_2, \ i_3=e_2+c_3f_2, \ j_3=f_2+f_3+c_3(-c_2+g_2)$.

  Things seem pretty bleak, but actually our standard form is close at
  hand. Let us return to $Q_2=XZ+ b_2Y^2+c_2YZ+d_2YW+e_2Z^2+f_2ZW$. We
  have $b_2=0$ since $p_2$ lies on $q_2$, so applying projective
  transformations which fix $p_1, \ p_2, \text{ and } Q_1$ , we can
  put $Q_2$ in the form $Q_2=XZ+YW$.  We can substitute in $d_2=1, \
  c_2=e_2=f_2=0$ in the equations above; the result is that we solve
  for $e_3, \ f_3, \ g_3 \ i_3 \text{ and } j_3$ in terms of $c_3$.
  Explicitly, we get $Q_3=WX-Y^2+c_3(WY+XZ)+h_3Z^2$.  Finally,
  applying projective transformations which fix $p_1, \ p_2 \text{ and }
  Q_1$ and map $Q_2$ to some quadric in the pencil $\left<Q_1,\
    Q_2 \right>$, we can put $Q_3$ in the form $Q_3=XW-Y^2$ or
  $Q_3=XW-Y^2+Z^2$, according as $h_3=0$ or not.  One can compute that
  if $Q_3=XW-Y^2$, the base locus of the net is not $0$-dimensional,
  so the standard form is as claimed.

\item $\{6,2\}$: In this case there are $2$ basepoints of the net, so
  without loss of generality we can assume these are $p_1=[1,0,0,0]$
  and $p_7=[0,1,0,0]$. There is a unique rank-$2$ quadric in the net,
  which we know is singular at $p_1$ and smooth at $p_7$. Therefore
  its equation is $Q_1=(a_1Y+b_1Z+c_1W)(d_1Z+e_1W)$, where the linear
  forms $a_1Y+b_1Z+c_1W$ and $d_1Z+e_1W$ are linearly independent. We
  can apply projective transformations fixing $p_1$ and $p_7$ to make
  this $Q_1=YZ$.

  Now $p_7$ is a multiple basepoint of the net, so there must be some
  quadric in the net singular there. It cannot be a double plane,
  since this would be singular at $p_1$ too, and by Assumption 1 only
  one quadric in our net may be singular at a given basepoint. So we
  can take $Q_2$ to be an irreducible reduced cone with vertex $p_7$.
  Such a cone is given by a form with no monomials involving $Y$: we
  can write it as $Q_2=a_2XZ+b_2XW+e_2Z^2+f_2ZW+g_2W^2$. Note that
  $a_2$, $b_2$ cannot both be zero, since $Q_2$ cannot be singular at
  $p_1$.

  Now, applying projective transformations fixing $p_1$, $p_7$ and
  $Q_1$, we can put $Q_2$ in the form $Q_{2a}=XW+Z^2$ or
  $Q_{2b}=XZ+W^2$, according to whether $b_2 \neq 0$ or $b_2=0$. (Note
  that no projective transformation fixing $p_1, \ p_7 \text{ and }
  Q_1$ takes the pencil $\left<Q_1,Q_{2a} \right>$ to the pencil
  $\left<Q_1, Q_{2b} \right>$.  For any such transformation would have
  to take $Q_{2a}$ to $Q_{2b}$ since they are the only quadrics in
  each pencil singular at $p_7$.  But it is easy to show that no such
  transformation fixing $p_1$, $p_7$ and $Q_1$ exists.)

  What of $Q_3$? We know it is a quadric containing $p_1$ and $p_7$,
  so the coefficients of $X^2 \text{ and } Y^2$ in $Q_3$ must be zero.
  Moreover, we can subtract a multiple of $Q_1$ to make the
  coefficient of $YZ$ equal to zero too. If $Q_2=Q_{2a}$, we can
  subtract a multiple of $Q_2$ to make the coefficient of $Z^2$ equal
  zero; if $Q_2=Q_{2b}$, we can arrange that the coefficient of $W^2$
  is zero. So we get two possibilities: $Q_{3a}=
  b_3XY+c_3XZ+d_3XW+g_3YW++i_3ZW+j_3W^2$, or $Q_{3b}=
  b_3XY+c_3XZ+d_3XW+g_3YW+ h_3Z^2+i_3ZW$.

  Our combinatorial classification showed that the curves $Q_2 \cap
  \{Y=0\}$ and $Q_3 \cap \{Y=0\}$ must have an intersection
  point of order 4 at $p_1$. As before, we can translate this
  condition into constraints on the coefficients of $Q_3$. In both
  cases, we get the conditions $c_3=d_3=i_3=0$.
      
  If $Q_2=Q_{2a}, \ Q_3 = Q_{3a}$, we $Q_{3a}= b_3XY+g_3YW+j_3W^2$. We
  can then apply projective transformations fixing $p_1$, $p_7$, $Q_1
  $ and $ Q_2$ to put it in the form $Q_{3a}=XY+W^2$. But now we note
  the following: the intersection $Q_{2a} \cap \{Z=0\}$ is a
  reducible conic $XW$, while $Q_{3a} \cap \{Z=0\}$ is a smooth conic
  whose tangent line at $p_7$ is $\{X=0\}$. So these $2$ curves have
  intersection multiplicity 3 at $p_7$, meaning that this net is
  actually of type $\{5,3\}$.

     It remains to consider $Q_2=Q_{2b}, \ Q_3=Q_{3b}$. In this case
     we get $Q_{3b}= b_3XY+g_3YW+h_3Z^2$, and admissible projective
     transformations put this in one of two forms: $Q_{3b}=XY+Z^2$ (if
     $g_3=0$) or $Q_{3b}^{\prime}=XY+YW+Z^2$ (if $g_3 \neq 0$). But in fact the
     resulting nets are projectively equivalent: the projective
     transformation $\phi \in PGL(4,k)$ with matrix

     \begin{align*}
       \phi &= \begin{pmatrix} 
         1&0&-\frac{1}{4}&1 \\
         0&1&0&0 \\
         0&0&1&0 \\
         0&0&-\frac{1}{2}&1
       \end{pmatrix}
     \end{align*}
     fixes $p_1$, $p_7$, $Q_1$ and $Q_2$, and one can write
     $Q_{3b}^{\prime}= \frac{1}{4}\phi(Q_1)+\phi(Q_{3b})$. So $\phi$ maps
     the net $\left<Q_1,Q_2,Q_{3b}\right>$ to the net
     $\left<Q_1,Q_2,Q_{3b}^{\prime}\right>$. Hence all extremal nets of type $\{6,2\}$
     have the standard form claimed.

   \item $\{5,3\}$: The argument in this case goes through exactly as
     in the previous one. We have $2$ basepoints, which we can choose
     to be $p_1=[1,0,0,0] $ and $p_6=[0,1,0,0]$; there is a unique rank-$2$
     quadric $Q_1$ in the net, which we can transform to $Q_1=YZ$;
     there is a unique quadric in the net $Q_2$ which is singular at
     $p_6$. Exactly the same argument as above shows that we can put
     this in the form $Q_2=XW+Z^2$, and then put $Q_3$ in the form
     $XY+W^2$. So this type has the standard form we claimed.

   \item $\{4,4\}_1$: In this case the net has a single rank-$2$
     quadric $Q_1$, which has multiplicity data $1^32^1 + 1^12^3$. We
     can apply projective transformations to put the basepoints at
     $p_1=[1,0,0,0] $ and $ p_5=[0,1,0,0]$, and to put $Q_1$ in the
     form $Q_1=ZW$.  Moreover, without loss of generality the
     plane $\{Z=0\}$ has the correct tangent direction at $p_1$, and
     $\{W=0\}$ has the correct tangent direction at $p_5$.

     Now take $2$ other quadrics $Q_2 $ and $ Q_3$ which, together
     with $Q_1$, span the net. We can write down quadratic forms
     defining these quadrics:

            \begin{align*}
              Q_2 &= a_2XY+b_2XZ+c_2XW+d_2YZ+e_2YW+f_2Z^2+g_2ZW+h_2W^2 \\
              Q_3 &= a_3XY+b_3XZ+c_3XW+d_3YZ+e_3YW+f_3Z^2+g_3ZW+h_3W^2 \\
            \end{align*}

            Since $Q_1$ is singular at both basepoints, the
            differentials $dQ_2$ and $dQ_3$ must be linearly
            independent at the basepoints, by Assumption 1.  In affine
            coordinates $\{X=1\}$ near $p_1$, their tangent spaces are
            $T_{p_1}Q_2 = \{a_2Y+b_2Z+c_2W=0\}, \ T_{p_1}Q_3 =
            \{a_3Y+b_3Z+c_3W=0\}$. Restricting to the plane $\{Z=0\}$
            these tangent spaces must coincide, which says that
            $a_2Y+c_2W $ and $ a_3Y+c_3W$ are linearly dependent. On
            the other hand, restricting to the plane $W=0$, the
            tangent spaces are transverse, implying that $a_2Y+b_2Z $
            and $ a_3Y+b_3Z$ are linearly independent. The analogous
            argument near $p_5$ says that $a_2X+d_2Z $ and $
            a_3X+d_3Z$ are linearly dependent, while $a_2X+e_2W $ and
            $ a_3X+e_3W$ are linearly independent. In particular we
            see that neither $a_2$ nor $a_3$ can be zero, so without
            loss of generality, we can divide through the two forms
            above to get $a_2=a_3=1$. Then linear dependence implies
            $c_2=c_3 $ and $ d_2=d_3$. Now scaling $Z $ and $ W$
            (which does not affect $Q_1$) we can assume $c_2=c_3=1, \
            d_2=d_3=1$.

            Now consider the intersection $Q_2 \cap Q_3 \cap \{Z=0\}$;
            this should consist of $p_1$ with multiplicity 3, and
            $p_5$ with multiplicity 1.  Setting $Z=0 $ and $ X=1$ in the
            forms defining $Q_2 $ and $ Q_3$, and setting
            $a_2=a_3=c_2=c_3=d_2=d_3=1$ as explained above, we get the
            forms

            \begin{align*}
              q_2 &= Y+W+e_2YW+h_2W^2 \\
              q_3 &= Y+W+e_3YW+h_3W^2 \\
            \end{align*}

            Setting $q_3=0$, again we can solve for $W$ as a power
            series in $Y$. Up to terms of order $4$, this is
            $W=Y(-1+(e_3-h_3)Y)$.  Substituting this in $q_2$, we get
            $q_2(Y)= Y^2 \left(-e_2+e_3+h_2-h_3\right)$, and this
            vanishes to order 3 at $Y=0$ if and only if
            $e_2-e_3=h_2-h_3$.

            Similarly consider the intersection $Q_2 \cap Q_3 \cap
            \{W=0\}$: this should consist of a simple point at $p_1$
            and a triple point at $p_5$. Setting $W=0 $ and $ Y=1$ in the
            forms defining $Q_2 $ and $ Q_3$, we get
            \begin{align*}
              q_2 &= X+Z+b_2XZ+f_2Z^2 \\
              q_3 &= X+Z+b_3XZ+f_3Z^2 \\
            \end{align*}

            By exactly the same reasoning as before, we get the
            equation $b_2-b_3=f_2-f_3$. So putting all these facts
            together, and subtracting multiples of $Q_1$ from $Q_2
            $ and $ Q_3$ to eliminate the monomials $ZW$ in each, we can
            write our quadrics as
            \begin{align*}
              Q_2 &= XY+b_2XZ+XW+YZ+e_2YW+f_2Z^2+h_2W^2 \\
              Q_3 &= XY+b_3XZ+XW+YZ+e_3YW+(f_2-b_2+b_3)Z^2+(h_2-e_2+e_3)W^2 \\
            \end{align*}

            But now
            \begin{align*}
              Q:= Q_2-Q_3 &=(b_2-b_3)XZ+(e_2-e_3)YW+(b_2-b_3)Z^2+(e_2-e_3)W^2\\
              &=(b_2-b_3)Z(X+Z) +(e_2-e_3)W(Y+W).
            \end{align*}

            Neither of the coefficients can be zero, since $Q$ is not
            reducible; scaling $X$ and $Z$ together and $Y$ and $W$
            together, we can put $Q$ in the form $Q=Z(X+Z)+W(Y+W)$,
            and changing coordinates $X \mapsto X+Z$, $Y \mapsto Y+W$,
            this becomes $Q=XZ+YW$. (Note that none of these changes
            affect $p_1$, $p_5$ or $Q_1$.) We will take $Q$ to be the
            second generator of our net, so we rename it $Q_2$: that
            is, we define $Q_2:=XZ+YW$. Note the following: the
            intersection $Q_1 \cap Q_2$ is a union of lines with total
            degree four: the line $\{Y=Z=0\}$, the line $\{X=W=0\}$,
            and the line $\{Z=W=0\}$ with multiplicity $2$. Any other
            quadric $Q$ which together with $Q_1$ and $Q_2$, spans the
            net, must pass through $p_1$ and $p_5$, hence must
            intersect $\{Z=W=0\}$ transversely at those two points. So
            the correct tangent direction at $p_1$ is the line
            $\{Y=Z=0\}$, and that at $p_5$ is the line $\{X=W=0\}$.

            Now suppose $Q_3$ is any other quadric which, together
            with $Q_1$ and $Q_2$, spans the net. Since it passes
            through $p_1$ and $p_5$, it has the form

            \begin{align*}
             Q_3 = a_3XY+b_3XZ+c_3XW+d_3YZ+e_3YW+f_3Z^2+g_3ZW+h_3W^2.
            \end{align*}

            We can subtract arbitrary multiples of $Q_1$ and $Q_2$
            without affecting anything, so we may assume that the
            coefficient $b_3$ of $XZ$ and the coefficient $g_3$ of
            $ZW$ both vanish. We know that if we restrict to the plane
            $Z=0$, the tangent line to $Q_3$ at $p_1$ should be the
            line $\{Y=0\}$. So $a_3 \neq 0$, $c_3=0$. Similarly
            restricting to $\{W=0\}$, the tangent line should be
            $\{X=0\}$, so $d_3=0$. So (dividing across by $a_3$) we
            get $Q_3 = XY+e_3YW+f_3Z^2+h_3W^2$. Neither of the
            coefficients $f_3$ or $h_3$ can be zero: if $f_3=0$, then
            $Q_1$, $Q_2$ and $Q_3$ all contain the line $\{X=W=0\}$;
            if $h_3=0$, they all contain $\{Y=Z=0\}$. By assumption
            our net has base locus of dimension 0, so this is
            forbidden. With this restriction, it is not difficult to
            see that for any values of the coefficients $e_3$, $f_3$,
            $h_3$, the net $\left<Q_1, Q_2, Q_3= XY+e_3YW+f_3Z^2+h_3W^2\right>$
            is projectively equivalent to the net $\left<Q_1, Q_2, Q_3=
            XY+Z^2+W^2\right>$, and this gives the standard form claimed.

        \item $\{4,4\}_2$: In this case the unique rank-$2$ quadric is
          smooth at both basepoints, so we can move the basepoints to
          $p_1=[1,0,0,0], \ p_5=[0,1,0,0]$, and transform the rank-$2$
          quadric to $Q_1=XY$. In this case the net contains a double
          plane $Q_2=L^2$ (where $L$ is a homogeneous linear form). It
          passes through both $[1,0,0,0]$ and $[0,1,0,0]$, so the
          coefficients of both $X$ and $Y$ in $L$ must vanish: hence
          by changing coordinates $Z \mapsto L(Z,W)$ (and $W \mapsto Z
          \text{ if } L=W$) we can assume that $Q_2=Z^2$.

          Now consider the third generator of the net, which by
          Assumption 1 must be a quadric smooth at both basepoints. We
          can write it as $Q_3=b_3XY+ c_3XZ
          +d_3XW+f_3YZ+g_3YW+h_3Z^2+i_3ZW+j_3W^2$ (the coefficients of
          $X^2 $ and $ Y^2$ are zero because $Q_3$ passes through $p_1
          $ and $ p_5$). Moreover, $Q_1 \cap Q_2 \cap Q_3$ should have
          $2$ points of multiplicity 4 at $p_1$ and $p_5$: this
          implies the double lines $Q_2 \cap \{Y=0\} $ and $ Q_2 \cap
          \{X=0\}$ should be tangent to the curves $Q_3 \cap \{Y=0\} $
          and $ Q_3 \cap \{X=0\}$ at $p_1, \ p_5$ respectively. In
          suitable affine coordinates near these points, the tangent
          lines to $Q_3$ inside these planes are defined by
          $c_3z+d_3w$ and $f_3z+g_3w$ respectively, so we get
          $d_3=g_3=0$ (and $c_3, \ f_3$ nonzero). Finally, replacing
          $Q_3$ by $Q_3-b_3Q_1-h_3Q_2$, we get $Q_3=c_3XZ
          +f_3YZ+i_3ZW+j_3W^2$.

          Note that of the four coefficients of $Q_3$, only $i_3$ can
          be zero: if $c_3$ were, $Q_3$ would be a cone with vertex
          $p_1$; if $f_3$ were, it would be a cone with vertex $p_5$;
          if $j_3$ were, $Q_3$ would be divisible by $Z$, which would
          make it a second rank-$2$ quadric in the net. If $i_3=0$, we
          get $Q_3=c_3XZ +f_3YZ+j_3W^2$: changing variables $X \mapsto
          c_3 X, \ Y \mapsto f_3Y, W \mapsto \sqrt{j_3} W$, we get
          $Q_3=XZ+YZ+W^2$ (without changing $p_1, \ p_5, \ Q_1 $ or $
          Q_2$). If $i_3$ is nonzero, we can do a similar rescaling of
          variables to make $c_3=f_3=i_3=j_3=1$, and $Q_3=XZ
          +YZ+ZW+W^2$. But then replacing $Q_3$ by $Q_3 + Q_2/4$, we
          get $Q_3=XZ+YZ+W^2+WZ+(Z/2)^2$, and finally changing
          variables $W \mapsto W+Z/2$, we get $Q_3=XZ+YZ+W^2$ again.
          So a net of this type containing a double plane can always
          be put in this standard form, as claimed.

 \item $\{4,4\}_3$: Again the unique rank-$2$ quadric in the
     net is smooth at both basepoints. We can move the basepoints by
     projective transformations to $p_1=[1,0,0,0], \ p_5=[0,1,0,0]$,
     and transform the rank-$2$ quadric to $Q_1=XY$.

     The net contains no double plane. Therefore the unique quadrics
     in the net singular at the basepoints $p_1$ and $p_5$ must be
     irreducible reduced cones with vertices at $p_1 $ and $ p_5$. Call
     these $Q_2$ and $Q_3$ respectively: then $C_2 := Q_2 \cap
     \{Y=0\}$ must be a reducible conic curve (that is, the union of
     $2$ lines in the plane $\{Y=0\}$, which may be equal), and
     similarly for $C_3:= Q_3 \cap \{X=0\}$. On the other hand
     $\Gamma_3:= Q_3 \cap \{X=0\} $ and $ \Gamma_2 := Q_2 \cap \{Y=0\}$
     are smooth conic curves in those planes, each meeting the
     reducible conic in the same plane in a single point of
     multiplicity 4. It follows that $C_2$ (resp.  $C_3$) is a double
     line, tangent at $p_1$ (resp. $p_5$) to the smooth conic
     $\Gamma_3$ (resp.  $\Gamma_2$).

     Let us write $Q_2=a_2YZ+b_2YW+c_2Z^2+d_2ZW+e_2W^2$,
     $Q_3=a_3XZ+b_3XW+c_3Z^2+d_3ZW+e_3W^2$. The restriction of $Q_2$
     (resp. $Q_3$) to $\{Y=0\}$ (resp. $\{X=0\}$) is a double line, so
     we get $d_2=\pm 2 \sqrt{c_2e_2},\ d_3= \pm 2 \sqrt{c_3 e_3}$.
     Rewriting, we have $Q_2=Y(a_2Z+b_2W) + (\gamma_2Z+\epsilon_2W)^2,
     \ Q_3 = X(a_3Z+b_3W) + (\gamma_3Z+\epsilon_3W)^2$, for some
     choice of square roots $\gamma_i, \ \epsilon_i$ of $ c_i, \ e_i \
     (i=1,\ 2)$. Now if the forms $\gamma_2Z+\epsilon_2W $ and $
     \gamma_3Z+\epsilon_3W$ are linearly dependent, then $Q_2$ and
     $Q_3$ would have an intersection point on the line $\{X=Y=0\}
     \subset Q_1$, which is impossible since the net has only $2$
     basepoints. Therefore they must be linearly independent, so we
     can change variables in $Z $ and $ W$ to make $Q_2=Y(a_2Z+b_2W)+Z^2,
     \ Q_3=X(a_3Z+b_3W)+W^2$. Now $Q_2 \cap \{X=0\}$ should be tangent
     to the double line $Q_3 \cap \{X=0\} = W^2$, so we get $a_2=0$;
     an identical argument gives $b_3=0$. Rescaling via $Y \mapsto
     b_2Y $ and $ X \mapsto a_3X$, we get $Q_2=YW+Z^2,\ Q_3=XZ+W^2$.
     Finally, we can swap $Q_2$ and $Q_3$, and our net has the
     standard form we claimed.    
     
   \item $\{4,2,2\}$: In this case we have $3$ distinct basepoints.
     By Lemma \ref{lemma_nondegeneracy} these do no lie on a line, so
     we can move them to $p_1=[1,0,0,0], \ p_5=[0,1,0,0], \
     p_7=[0,0,1,0]$. The combinatorial classification shows that $Q_1$
     can be taken to be a rank-$2$ quadric $P_1 \cup P_2$, where $P_1$
     is a plane passing through $p_1$, but not through $p_5$ or $p_7$,
     and $P_2$ is a plane passing through $p_5$ and $p_7$, but not
     $p_1$. So we can write these as $P_1=b_1Y+c_1Z+d_1W, \
     P_2=a_2X+d_2W$, with $b_1$, $c_1$, $a_2 \neq 0$. Changing
     coordinates $X \mapsto a_2X+d_2W$, $Y \mapsto b_1Y+d_1W$, $Z
     \mapsto c_1Z$ (which does not affect $p_1$, $p_5$ or $p_7$) we
     obtain $Q_1=X(Y+Z)$.
        
        Now for $Q_2$. It is a rank-$2$ quadric, consisting of a plane
        $\Pi_1$ passing through $p_1 $ and $ p_5$, and a plane
        $\Pi_2$ passing through $p_1 $ and $ p_7$. So we have
        $\Pi_1 = c_1Z+d_1W, \ \Pi_2=b_2Y+d_2W$, with $c_1 $ and $ b_2$
        nonzero; dividing out, we can assume these coefficients both
        equal $1$.  Now, each of these $2$ planes should contain the
        tangent line at $p_1$ which is the first basepoint infinitely
        near to $p_1$; hence, in terms of embedded tangent spaces,
        that tangent line is the intersection $\Pi_1 \cap \Pi_2$.
        Moroever, we know that the plane $P_1$ defined above must also
        contain that tangent line. That means the lines $P_1 \cap
        \Pi_1 = \{-Y+d_1W=0\} $ and $ P_1 \cap \Pi_2 = \{Y+d_2W=0\}$
        are equal, hence $d_2=-d_1$. Now applying the transformation
        $Y \mapsto Y-d_1W$, $Z \mapsto Z+d_1W$ we get $Q_2=YZ$, and
        $p_1, \ p_5, \ p_7 $ and $ Q_1$ are unchanged.

        Finally we must deal with $Q_3$. We know it passes through
        $p_1$, $p_5$, and $p_7$, so the coefficients of $X^2, \ Y^2,
        $ and $ Z^2$ must be zero. So write
        $Q_3=a_3XY+b_3XZ+c_3XW+d_3YZ+e_3YW+f_3ZW+g_3W^2$. Moreover, we
        know the tangent direction $Q_3$ must have at the $3$
        basepoints. At $p_1$, the correct tangent line is that shared
        by $\Pi_1 $ and $ \Pi_2$ above, namely $\{Y=Z=0\}$.  Setting
        $X=1$ in the equation of $Q_3$, we get $a_3Y+b_3Z+c_3W +
        \text{ (quadratic terms)}$. So we get the condition $c_3=0$.
        Now consider $p_5$: the correct tangent direction there is
        that shared by the planes $P_1$ and $\Pi_1$, and that is
        $\{X=Z=0\}$. Setting $Y=1$ in the equation of $Q_3$, we get
        $a_3X+d_3Z+e_3W + \text{ (quadratic terms)}$, so the condition
        we get is $e_3=0$.  Finally looking at $p_7$, the correct
        tangent direction is that shared by $P_1$ and $\Pi_2$, and the
        same argument gives the condition $f_3=0$. So these three
        conditions give us $Q_3=a_3XY+b_3XZ+d_3YZ+g_3W^2$. But now
        replacing $Q_3$ by $Q_3-d_3Q_2-a_3Q_1$, we can eliminate the
        monomials $YZ$ and $XY$, giving $Q_3=b_3XZ+g_3W^2$. Neither
        coefficient can be zero --- if $b_3$ were zero, $Q_3$ would be
        a double plane, hence singular at $p_1$, but this would
        violate Assumption 1 since $Q_1$ is singular there; if $g_3$
        were zero, then $Q_3$ would be a third rank-$2$
        quadric in the net. So both are nonzero; dividing across by
        $b_3$ and scaling $W$ (which does not affect the basepoints or
        $Q_1 $ and $ Q_2$) we get $Q_3=XZ+W^2$, as claimed.

        \item $\{3,3,2\}_1$: Again we can put the $3$ basepoints at
          $p_1=[1,0,0,0], \ p_4=[0,1,0,0], \ p_7=[0,0,1,0]$. In this
          case, the rank-$2$ quadrics in the net have multiplicity
          data $Q_1=1^33^1+2^33^1$, $Q_2=1^22^2+1^12^13^2$. So they
          have equations $Q_1=(b_1Y+d_1W)(a_2X+d_2W)$,
          $Q_2=(\gamma_1Z+\delta_1W)W$.  None of the coefficients
          $b_1, \ a_2, \ \gamma_1$ can be zero, otherwise the
          corresponding planes would pass through more basepoints than
          specified by the combinatorial classification. So by
          changing coordinates $X \mapsto a_2X+d_2W$, $Y \mapsto b_1Y+d_1W$,
          $Z'=\gamma_1Z+\delta_1W)$, we obtain $Q_1=XY$, $Q_2=ZW$.

          Now consider $Q_3$, any quadric in the net which forms a
          basis together with $Q_1$ and $Q_2$. Such a $Q_3$ must pass
          through $p_1, \ p_4 $ and $ p_7$. Moreover, $Q_1$ is singular
          at one $\P^3$-basepoint, and $Q_2$ is singular at the other
          $2$, so $Q_3$ is smooth at the base locus, and has the
          correct tangent direction at each. But $Q_1$ and $Q_2$
          define the correct tangent direction at $p_1$ and $p_4$.
          Applying these conditions to the quadratic form defining
          $Q_3$, we see that the coefficients of the monomials $X^2$,
          $Y^2$, $Z^2$, $XW $ and $ YW$ must all be zero. So we can write
          $Q_3=a_3XY+b_3XZ+c_3YZ+d_3ZW+e_3W^2$.

          Now the above facts (about smoothness of $Q_3$ at the base
          locus, and its tangent directions there) hold for any
          quadric in the net outside the pencil spanned by $Q_1$ and
          $Q_2$. In particular they remain true if we replace $Q_3$ by
          $Q_3-a_3Q_1-d_3Q_2$. So without loss of generality we obtain
          $Q_3=b_3XZ+c_3YZ+e_3W^2$. Now we see that $e_3$ must be nonzero, for
          otherwise $Q_3$ would be reducible. Also, in affine
          coordinates near $p_1$ and $p_4$ the tangent spaces to $Q_3$
          are given by $b_3z=0$ and $c_3z=0$ respectively. Smoothness at
          these points tells us that $b_3$ and $c_3$ are nonzero. So all
          three coefficients are nonzero; scaling the coordinates we
          get $Q_3=XZ+YZ+W^2$, as claimed.

 \item $\{3,3,2\}_2$: The combinatorial classification tells us
        in this case that one of the rank-$2$ quadrics in the net ---
        let us call it $Q_1$ --- is the union of a plane $P_1$ passing
        through $p_1 $ and $ p_4$, and a plane $P_2$ passing through
        $p_1 $ and $p_7$.  These are given by forms $P_1=c_1Z+d_1W, \
        P_2=b_2Y+d_2W$, and exactly as in the previous case we can
        transform these to $P_1=Z, \ P_2 = Y$.  So $Q_1=YZ$. Similarly
        the other rank-$2$ quadric in the net --- call it $Q_2$ --- is the
        union of a plane $\Pi_1$ through $p_1 $ and $ p_4$, and a
        plane $\Pi_2$ through $p_4 $ and $ p_7$; by exactly the same
        argument, we can put this in the form $Q_2=X(Z+W)$.

          What of $Q_3$? As in the previous case, we know the
          coefficients of the monomials $X^2, \ Y^2 $ and $ Z^2$ in $Q_3$
          must be zero. Also, just as before, we can compute the
          shared tangent directions of components of $Q_1 $ and $ Q_2$ at
          the basepoints: this tells us that the coefficients of $YW
          $ and $ ZW$ in $Q_3$ are zero, and those of $XZ $ and $ XW$ must
          be equal. So we get $Q_3=a_3XY+b_3(XZ+XW)+c_3YZ+d_3W^2$. But
          now replacing $Q_3$ by $Q_3-b_3Q_2-c_3Q_1$, we get
          $Q_3=a_3XY+d_3W^2$.  Just as in the previous case, neither
          coefficient can be zero, so we can rescale via $X \mapsto
          a_3X, (W,Z) \mapsto \sqrt{d_3}(W,Z)$ (without moving the
          basepoints or $Q_1 $ and $ Q_2$) to get $Q_3=XY+W^2$, as
          claimed.
          
        \item $\{2,2,2,2\}$: First note that the 4 $\P^3$-basepoints
          of the net cannot be coplanar. The proper transform of such
          a plane would have class $h_{1357}$, but the combinatorial
          classification shows there is an effective class $h_{1257}$;
          the corresponding planes in $\P^3$ must then be equal, which
          means that in fact the class $h-e_1-e_2-e_3-e_5-e_7$ would
          be effective, which is impossible.

          We know also that no 3 of the $\P^3$-basepoints are
          collinear. So we can move them to the coordinate points of
          $\P^3$: $p_1=[1,0,0,0]$, $p_3=[0,1,0,0]$, $p_5=[0,0,1,0]$,
          $p_7=[0,0,0,1]$. The combinatorial classification shows that
          the multiplicity data of the rank-$2$ quadrics in the net
          are as follows: $Q_1=1^12^13^2+1^12^14^2$,
          $Q_2=1^23^14^1+2^23^14^1$, $Q_3=1^22^2+3^24^2$. But then the
          components of $Q_1 $ and $ Q_2$ are determined: we have
          $Q_1=ZW$, $Q_2=XY$. We also get $Q_3=(aX+bY)(cZ+dW)$, with
          $a, \ b, \ c, d$ all nonzero. But then we can scale the
          coordinates (without changing the $p_i$ or $Q_1 $ and $
          Q_2$) to get $Q_3=(X+Y)(Z+W)$, as claimed.

        \item $\{1,1,1,1,1,1,1,1\}$: The combinatorial classification
          from the last section showed that the four points
          $\{p_1,p_2,p_3,p_5\}$ do not lie in a plane in $\P^3$, so we
          can move them to the coordinate points: $p_1=[1,0,0,0]$,
          $p_2=[0,1,0,0]$, $p_3=[0,0,1,0]$, $p_5=[0,0,0,1]$. We know
          that $p_4$ (resp. $p_6$, $p_7$) lies in the plane spanned by
          $\{p_1,p_2,p_3\}$ (resp.  $\{p_1,p_2,p_5\}$,
          $\{p_1,p_3,p_5\}$) so we have $p_4=[x_4,y_4,z_4,0]$,
          $p_6=[x_6,y_6,0,w_6]$, $p_7=[x_7,0,z_7,w_7]$, and the
          coordinates $x_i$, $y_j$, $z_k$, $w_l$ are all nonzero
          (since otherwise we would get $3$ collinear basepoints,
          which is forbidden). Normalising, we can write
          $p_4=[1,y_4,z_4,0]$, $p_6=[1,y_6,0,w_6]$,
          $p_7=[1,0,z_7,w_7]$. 

          What of $p_8$? We know it does not belong to any of the
          planes $\{Y=0\}$, $\{Z=0\}$, $\{W=0\}$, since each of those
          already contains $4$ basepoints. So it has coordinates
          $p_4=[x_8,y_8,z_8,w_8]$, with $y_8z_8w_8 \neq 0$. On the
          other hand, we know that $p_8$ lies in the plane spanned by
          $\{p_2,p_3,p_5\}$, so it must have $x_8=0$. Applying the
          projective transformation $[X,Y,Z,W] \mapsto [X,
          Y/y_8,Z/z_8,W/w_8]$ to $\P^3$, we bring $p_8$ to
          $[0,1,1,1]$, without moving $p_1$, $p_2$, $p_3$, $p_5$, or
          changing the form of $p_4$, $p_6$, $p_7$ above.

          We know from the combinatorial classification that the points
          $\{p_1,p_4,p_5,p_8\}$ are coplanar. This is equivalent to
          the determinant of the matrix 
          \begin{center}
            \[ \left( \begin{array}{cccc}
                1 & 0 & 0 & 0 \\
                1 & y_4 & z_4 & 0 \\
                0 & 0 & 0 & 1 \\
                0 & 1 & 1 & 1 \end{array} \right)\]
          \end{center}

          (whose rows are the homogeneous coordinates of the four
          points) vanishing, which holds if and only if $y_4=z_4$.
          Similar arguments show we must have $y_6=w_6$ and $z_7=w_7$.

          Next, we use the fact that the points $\{p_1,p_4,p_6,p_7\}$
          are coplanar. That means the determinant of the corresponding
          matrix must vanish: this determinant is $-2y_4y_6z_7$, and
          we know $y_4$, $y_6$, $z_7$ are all nonzero. This shows that
          an extremal net of this type can only exist if the
          characteristic of the base field is $2$.

          To find the standard form in the case of characteristic $2$,
          we now use that the points $\{p_5,p_6,p_7,p_8\}$ are
          coplanar. Again we use vanishing of the determinant of the
          corresponding matrix: this determinant is $y_6+z_7$, so we
          get $y_6=z_7$. A similar argument shows that $y_4=y_6$. So
          our points have coordinates $p_4=[1,\xi,\xi,0]$,
          $p_6=[1,\xi,0,\xi]$, $p_7=[1,0,\xi,\xi]$, for some nonzero
          $\xi \in k$. Applying the projective transformation
          $[X,Y,Z,W] \mapsto [X,Y/\xi,Z/\xi,W/\xi]$, the points $p_4$,
          $p_6$, $p_7$ are transformed to $p_4 = [1,1,1,0]$,
          $p_6=[1,1,0,1]$, $p_7=[1,0,1,1]$, and the other five points
          are left fixed.

          Finally, consider the equations of the planes containing $4$
          of the basepoints. The plane containing
          $\{p_1,p_2,p_3,p_4\}$ has equation $W=0$, and the plane
          containing $\{p_5,p_6,p_7,p_8\}$ has equation $X+Y+Z=0$.
          This gives a rank-$2$ quadric $Q_1=(X+Y+Z)W=0$ in the net.
          The plane containing $\{p_1,p_2,p_5,p_6\}$ has equation
          $Z=0$, and the plane containing $\{p_3,p_4,p_7,p_8\}$ has
          equation $X+Y+W=0$, giving a rank-$2$ quadric
          $Q_2=(X+Y+W)Z=0$ in the net. The plane containing
          $\{p_1,p_3,p_5,p_7\}$ has equation $Y=0$, and the plane
          containing $\{p_2,p_4,p_6,p_8\}$ has equation $X+Z+W=0$,
          giving a rank-$2$ quadric $Q_3=(X+Z+W)Y=0$ in the net. This
          gives the standard form claimed.
    \end{enumerate}

\section{Extremal fibrations and extremal quartics} \label{section_extremalquartics}

In this section we assume that the characteristic of the ground field
$k$ is not $2$. (In particular, our remarks do not apply to the
extremal net of type $\{1,1,1,1,1,1,1,1\}$.) Suppose we are given a
net $N$ of quadrics in $\P^3$ with some fixed basis, say $N= \left<
\lambda_1Q_1+\lambda_2Q_2+\lambda_3Q_3 \right>$. The {\it discriminant
  form} $\Delta_N = \text{ det} \left(
\lambda_1Q_1+\lambda_2Q_2+\lambda_3Q_3 \right)$ defines a quartic curve
in the plane $N \iso \P^2$. It seems reasonable to expect that
extremality of the net $N$ in the sense used heretofar should
correspond to some extremality property of the quartic $N$.

To explain the correspondence, we first note that there is a natural
connection between plane quartic curves and the root system $E_7$.  To
an isolated hypersurface singularity one can associate in a natural way
a root system (see \cite[Chapter 4]{Arnold1988} for details). For
plane quartics, the ranks of the root systems associated to its
various singular points sum to at most $7$, and in this case the direct
sum of the root systems is a rank-$7$ root subsystem of $E_7$. So one
can hope that for an extremal net $N$ the quartic $\Delta_N$ is
extremal, in the sense that the associated root system has rank $7$.
Indeed, it seems natural to expect in this case that the root system
associated to $N$ in Table \ref{table_intro} and that associated to
$\Delta_N$ should in fact be the same. This is what we verify below.

The following table gives, for each type of extremal net $N$, a
defining equation for its discriminant quartic $\Delta_N$ and the root
system associated to the singularities of $\Delta_N$. (See for
instance \cite{bruce1980} for details on how to identify root systems
of singularities from equations.) Here $\lambda_1$, $\lambda_2$,
$\lambda_3$ are homogeneous coordinates on the net $N \iso \P^2$.

\begin{center}
\begin{tabular}{lll}
{\bf Type} & {\bf $\Delta_N$} & Singularities of $\Delta_N$\\ 
\toprule 
$\{8\}_1$ &
$\lambda_2 \left(
4\lambda_1\lambda_2^2+\lambda_2\lambda_3^2+4\lambda_3^3 \right)$
& $E_7$\\ 
$\{8\}_2$ &
$\lambda_2^4+2\lambda_1\lambda_2^2\lambda_3+\lambda_1^2\lambda_3^2+4\lambda_3^4$
& $A_7$\\ 
$\{4,4\}_1$ &
$\left(\lambda_2^2-\lambda_1\lambda_3+2\lambda_3^2\right) \left(\lambda_2^2-\lambda_1\lambda_3-2\lambda_3^2\right) $
& $A_7$\\ 
$\{6,2\}$ & $\lambda_2 \lambda_3 \left( \lambda_1\lambda_2-\lambda_3^2
\right)$& $D_6 + A_1$\\
$\{4,4\}_2$ & $\lambda_1\lambda_3 \left(\lambda_1\lambda_2-\lambda_3^2
\right)$
& $D_6+A_1$\\ 
$\{5,3\}$ & $\lambda_2 \left( \lambda_1^2\lambda_2-4\lambda_3^3 \right)$
& $A_5+A_2$\\ 
$\{3,3,2\}_1$ & $\lambda_1\left(\lambda_1\lambda_2^2-4\lambda_3^3 \right)$
& $A_5+A_2$\\
 $\{4,2,2\}$
& $\lambda_1\lambda_2\lambda_3 \left(\lambda_1+\lambda_2\right)$
& $D_4+3A_1$\\ 
$\{4,4\}_3$ &
$\lambda_2\lambda_3 \left( \lambda_1^2-\lambda_2\lambda_3 \right)$ & $2A_3+A_1$\\
$\{3,3,2\}_2$ &
$\lambda_1\lambda_2 \left(\lambda_1\lambda_2+4\lambda_3^2\right)$
&$2A_3+A_1$\\ 
$\{2,2,2,2\}$ &
$\lambda_1\lambda_2 \left(\lambda_1\lambda_2-4\lambda_3^2\right)$ &$2A_3+A_1$\\
\end{tabular}
\end{center}

We observe that in each case the root system associated to $\Delta_N$
is the same as that associated to $N$ in Table \ref{table_intro}. It
would be interesting to find an explanation for this correspondence.

\renewcommand{\baselinestretch}{1.62}

\small
\sc DPMMS, Wilberforce Road, Cambridge CB3 0WB, United Kingdom

Leibniz Universit\"at Hannover, Institut f\"ur Algebraische Geometrie,
Welfengarten 1, D-30167, Germany

{\it Email address:} {\tt artie@math.uni-hannover.de}


\begin{thebibliography}{99}

\bibitem{Arnold1988} V. Arnold, S. Guse{\u\i}n-Zade,
  A. Varchenko. {\it Singularities of differentiable
    maps. {V}ol. {II}.} Birkh\"auser (1988).

\bibitem{Barthetal2003} W. Barth, K. Hulek, C. Peters, A Van De
  Ven. {\it Compact Complex Surfaces. } $2^\text{nd}$ ed. Springer
  (2003). 

\bibitem{bruce1980} J. Bruce, P. Giblin. A stratification of the space
  of plane quartic curves. Proc. London Math. Soc. (3) {\bf 42}
  (1981), 270--298.


\bibitem{CossecDolgachev1989} F.Cossec, I. Dolgachev. {\it Enriques
    Surfaces I. } Birkh\"auser (1989).

\bibitem{DolgachevOrtland1988} I. Dolgachev, D. Ortland. Point sets in
  projective spaces and theta functions. {\it Ast\'erisque }{\bf 165} (1988).

\bibitem{DolgachevNotes} I. Dolgachev. Topics in Classical Algebraic
  Geometry. Part I. Available from {\tt
    http://www.math.lsa.umich.edu/\~{}idolga/lecturenotes.html}


\bibitem{GorbOnisVin1994} V. Gorbatsevich, A. Onischik,
  E. Vinberg. {\it Lie Groups and Lie Algebras III.} Encyclopaedia of
  Mathematical Sciences 41. Springer (1994).

\bibitem{Hartshorne1977} R. Hartshorne. {\it Algebraic Geometry. }
  Springer (1977).

\bibitem{HulekKloosterman2008} K. Hulek, R. Kloosterman. Calculating
the Mordell--Weil rank of elliptic threefolds and the cohomology of
singular hypersurfaces. arXiv: 0806.2025.

\bibitem{Humphreys1992} J. Humphreys. {\it Reflection Groups and
  Coxeter Groups. } Cambridge (1992).


\bibitem{Kodaira1963} K. Kodaira. On compact analytic surfaces
  II, III. {\it Ann. of Math. }{\bf 77} (1963), 563--626; {\bf 78}
  (1963), 1--40.


\bibitem{Lang1991} W. Lang, Extremal rational elliptic surfaces in
  characteristic $p$. I. Beauville surfaces. {\it Math. Z. } {\bf 207}
  (1991), 429--437.

\bibitem{Lang1994} W. Lang, Extremal rational elliptic surfaces in
  characteristic $p$. II. Surfaces with three or fewer singular
  fibres. {\it Ark. Math. } {\bf 32} (1994), 423--448.


\bibitem{MirandaPersson1986} R. Miranda and U. Persson, On extremal
  rational elliptic surfaces. {\it Math. Z. } {\bf 193} (1986),
  537--558.


\bibitem{Mukai2001} S. Mukai. Counterexample to Hilbert's fourteenth
  problem for the 3-dimensional additive group. RIMS Kyoto preprint
  \#1343.

 \bibitem{PrendergastSmith2009} A. Prendergast-Smith. Extremal
   rational elliptic threefolds. Ph.D. thesis, University of
   Cambridge, 2009.

 \bibitem{Totaro2008} B. Totaro. Hilbert's fourteenth problem over
   finite fields, and a conjecture on the cone of curves. {\it
     Compos. Math. } {\bf 144} (2008), 1176--1198.

\end{thebibliography}
\end{document}